\documentclass[12 pt]{article}

\usepackage{helvet} 

\usepackage{graphicx}  
\usepackage{amssymb}
\usepackage{amsmath}
\usepackage{subcaption}
\usepackage{wrapfig}
\usepackage{multirow}

\usepackage{rotating}
\usepackage{tikz}
\usetikzlibrary{arrows,
                positioning,
                shapes.geometric}

\usepackage{adjustbox}
                
\usepackage{hyperref}
\hypersetup{
    colorlinks=true,
    linkcolor=blue,
    citecolor=blue,
    filecolor= magenta,      
    urlcolor= magenta 
}

\newcommand{\Phit}{\tilde{\Phi}}
\newcommand{\es}{\epsilon_s}
\newcommand{\ep}{\epsilon_p}
\newcommand{\emm}{\epsilon_m}
\newcommand{\ez}{\epsilon_0}
\newcommand{\rr}{{\mathbf r}}
\newcommand{\s}{{\mathbf s}}
\newcommand{\nn}{{\mathbf n}}

\newcommand{\R}{\mathbb{R}}

\usepackage[margin=1.0in]{geometry}

\title{An Efficient Finite Element Iterative Method \\ for Solving a Nonuniform Size Modified Poisson-Boltzmann Ion Channel Model }
\author{Dexuan Xie\thanks{Department of Mathematical Sciences, University of Wisconsin-Milwaukee, Milwaukee, WI, 53201-0413, USA}}
\date{}	

\begin{document}


\maketitle

\begin{abstract}
In this paper, a nonuniform size modified  Poisson-Boltzmann ion channel (nuSMPBIC) model is presented as a nonlinear system of an electrostatic potential and multiple ionic concentrations. It mixes nonlinear algebraic equations with a Poisson boundary value problem involving Dirichlet-Neumann mixed boundary value conditions and a membrane surface charge density to reflect the effects of ion sizes and membrane charges on electrostatics and ionic concentrations. To  overcome the difficulties of strong singularities and exponential nonlinearities, it is split into three submodels with a solution of Model 1 collecting all the singular points and Models 2 and 3 much easier to solve numerically  than the original nuSMPBIC model. A damped two-block iterative method is then presented to solve Model 3, along with a novel modified Newton iterative scheme for solving each related nonlinear algebraic system. To this end, an effective nuSMPBIC finite element solver is derived and then implemented as a program package that works for an ion channel protein with a three-dimensional molecular structure and a mixture solution of multiple ionic species. Numerical results for a voltage-dependent anion channel (VDAC) in a mixture of four ionic species demonstrate a fast convergence rate of the damped two-block iterative method, the high performance of the  software package, and the importance of considering nonuniform ion sizes. Moreover, the nuSMPBIC model is validated by the anion selectivity property of VDAC. 
\end{abstract}


\section{Introduction}

The Poisson-Boltzmann equation (PBE) is one widely-used  dielectric continuum model for the calculation of electrostatic solvation free energies \cite{chen2011mibpb,honig95,jurrus2018improvements,luo2002accelerated,roux99,xiePBE2013}. But it cannot distinguish the two ions with the same charge, such as cations Na$^+$ and K$^+$, since it treats ions as volume-less points. Thus, it may work poorly in  the applications in which ion sizes have impact on electrostatics and ionic concentrations, especially in the simulation of ion transport across a membrane  via an ion channel pore. Thus, it is important to develop size modified PBE (SMPBE)  models.

The first SMPBE  model was reported in 1997 for an asymmetric electrolyte consisting of two ionic species under the assumption that all the ions have the same size \cite{borukhov1997steric}. Since then, several SMPBE models had been developed for the case of a protein in an ionic solvent. For example,  one SMPBE model for a protein in a solution of multiple ionic species was reported in \cite{BoLi2009a}. Since it treats all the ions and water molecules as cubes, this model has a void problem (i.e.,  there exist cavities among ions and water molecules) in the nonuniform ion size case, breaking down the required size constraint conditions. This  drawback was fixed by introducing a concentration of voids model in \cite[Eq. (10)]{liu2014poisson}, from which another SMPBE model could be constructed. But how to estimate a concentration of voids is still a puzzle since the voids can have different shapes and different volumes.  Besides, these two SMPBE models suffer a redundancy problem because they use a  concentration of water molecules to describe the water solution that has been treated as a continuum dielectric in their constructions. To fix these  drawbacks, several improved SMPBE models were reported in \cite{Li_Xie2014b, li2019analysis, nuSMPBE2017,Ying-Xie2015}, along with their effective finite element solvers and a web server  \cite{SMPBS2016}.  Most recently, a size modified Poisson-Boltzmann ion channel (SMPBIC) model and its effective finite element solver were reported in \cite{SMPBEic2019}. While it can partially reflect the effect of distinct ion sizes due to setting each ion to have an average of ion sizes, this SMPBIC model still cannot distinguish the two ions with the same charge. As the continuation of this previous work, the purpose of this work is to develop a nonuniform SMPBIC (nuSMPBIC) model and its finite element solver for an ion channel protein in a mixture of multiple ionic species. 

However, the nuSMPBIC model is much more difficult to solve numerically than the SMPBIC model since it is a nonlinear system  mixing $n$ nonlinear algebraic equations with one Poisson dielectric interface boundary value problem and involves two physical domains --- a simulation box domain for potential functions and a solvent domain for ionic concentration functions, not to mention its stronger singularities and stronger nonlinearities than the case of a protein surrounded by an ionic solvent.  Here $n$ is the number of ionic species; the algebraic equations describe ion size constraint conditions;  and the Poisson problem involves Dirichlet-Neumann mixed boundary value conditions and a membrane surface charge density to reflect membrane charge effects.  A  nuSMPBIC solution gives $n$ ionic concentration functions $c_i$ and an electrostatic potential function $u$. 
Since $c_i$ and $u$  are defined in two different domains, they belong to two different finite element function spaces, producing a two-domain issue  that we must deal with during the development of a nuSMPBIC finite element solver. Currently, this issue was simply treated by letting $c_i$ belong to the function space of $u$ through setting $c_i$ to be zero at the mesh points outside a solvent region  \cite{nuSMPBE2017}. But this simple treatment may cause a large error disturbance since it produces an artificial boundary layer around the solvent region, where $c_i$ may have large values due to strong atomic charges on a part of solvent domain boundary --- a part of an ion charge molecular surface including the mostly charged charge ion channel pore surface. To avoid such a boundary layer error, in this work, we treat this two physical domain issue directly through  constructing two finite element function spaces: one for concentration functions and the other for potential functions. We then construct two communication operators to directly carry out operations involving  ionic concentrations and potential functions. 

Following what was done in \cite{SMPBEic2019}, we will overcome the difficulties of solution singularities through splitting the nuSMPBIC model into three submodels, called Models 1, 2 and 3. While Models 1 and 2 are the same as those reported in  \cite{SMPBEic2019}, Model 3 is a nonlinear system mixing the $n$ nonlinear algebraic equations with an interface boundary value problem. A solution of Model 1 gives a potential component, $G$, induced by the atomic charges of an ion channel protein, in an algebraic expression that collects all the solution singularity points while a solution of Model 2 gives a potential component, $\Psi$, induced by membrane charges and Dirichlet boundary and interfaces values.  Solving Model 3 gives $n$ ionic concentration functions $c_i$ and a potential component function, $\Phit$, induced by ionic charges from the solvent region $D_s$. Model 2 has been solved efficiently in \cite{SMPBEic2019}. Hence, we only need to develop a  nonlinear iterative scheme for solving Model 3 in this work.

The classical nonlinear successive over-relaxation (SOR) iterative technique \cite{OR1970} is often used to develop iterative schemes for solving a system of nonlinear equations because it can separate each equation from the others so that each equation can be solved one-by-one to achieve the goal of reducing computational complexity and computer memory requirement. It was applied to the construction of a nonlinear SOR-Newton iterative scheme for solving a nonuniform SMPBE model in the case of a protein surrounded by an ionic solvent  \cite{nuSMPBE2017}.  Thus, we started with an adoption of this scheme as a Model 3 solver. Unfortunately,  this adapted SOR-Newton scheme was found numerically not to work well due to a slow rate of convergence because an ion channel protein has a much more complicated  geometry and much stronger atomic charges than the protein case considered in \cite{nuSMPBE2017}. Hence, developing efficient iterative schemes for solving Model 3 becomes one key step in the construction of an effective nuSMPBIC finite element solver. 

During the search for an efficient Model 3 solver, we discovered that under a linear finite element  framework, all the nonlinear algebraic equations of Model 3 can be split from a large set of $n N_h $ nonlinear  algebraic equations  into $N_h$ small sets with each set containing only $n$ nonlinear algebraic equations. Here $N_h$ is the number of mesh points of  a solvent region mesh. Since $n$ is typically very small (such as 2, 3, or 4) in practice, each small set can be quickly solved as a small nonlinear algebraic system. This discovery motivated us to divide  the equations of Model 3 into two blocks --- Block 1 consists of all the nonlinear algebraic equations and Block 2 contains the linear boundary problem only. We then construct a novel modified Newton iterative scheme to solve the nonlinear algebraic equations of Block 1 quickly, resulting in an efficient damped two-block iterative method for solving Model 3. In addition, we construct a linearized SMPBIC model and use its finite element solution as a good initial iterate of the  damped two-block iterative method. Consequently, the construction of an efficient nuSMPBIC finite element solver is completed.  
 
Finally, we implemented this new nuSMPBIC finite element solver in Python and Fortran as a software package based on  the state-of-the-art finite element library from the FEniCS project \cite{fenics-book} and the SMPBIC program package \cite{SMPBEic2019}. To demonstrate the performance  of the new nuSMPBIC software package, we did numerical tests using  a crystallographic three-dimensional molecular structure of a murine voltage-dependent anion channel 1 (mVDAC1)  \cite{ujwal2008crystal} in a mixture solution of four ionic species. 
Note that this mVDAC1 protein (PBD ID: 3EMN) is known to be in the open state conformation with an anion selectivity property. Hence, it can serve as a good test case for assessing our nuSMPBIC model.  Numerical results demonstrate a fast convergence rate of our damped two-block iterative method, the high performance of our  nuSMPBIC software package, and the importance of considering nonuniform ion sizes. They also show that  the membrane surface charges have impact on the electrostatic potential and ionic concentrations. Moreover, our SMPBIC model was validated by the mVDAC1 anion selectivity property. 
 
The remaining part of the paper is arranged as follows. In Section 2, we present the  nuSMPBIC model. In Section~3, we split the nuSMPBIC model into three submodels. In Section~4,  we present the damped two-block iterative method.  In Section 5, we present  the Newton iterative scheme for solving each nonlinear algebraic system arisen from the damped two-block iterative method. In Section 6, we construct a good initial iterate for the damped two-block iterative method. In Section 7,  we report the nuSMPBE program package  and  numerical test results. Finally, conclusions are made in Section 8. 

\begin{figure}
  \begin{center}
  \vspace{-15mm}
    \includegraphics[width=0.4\textwidth]{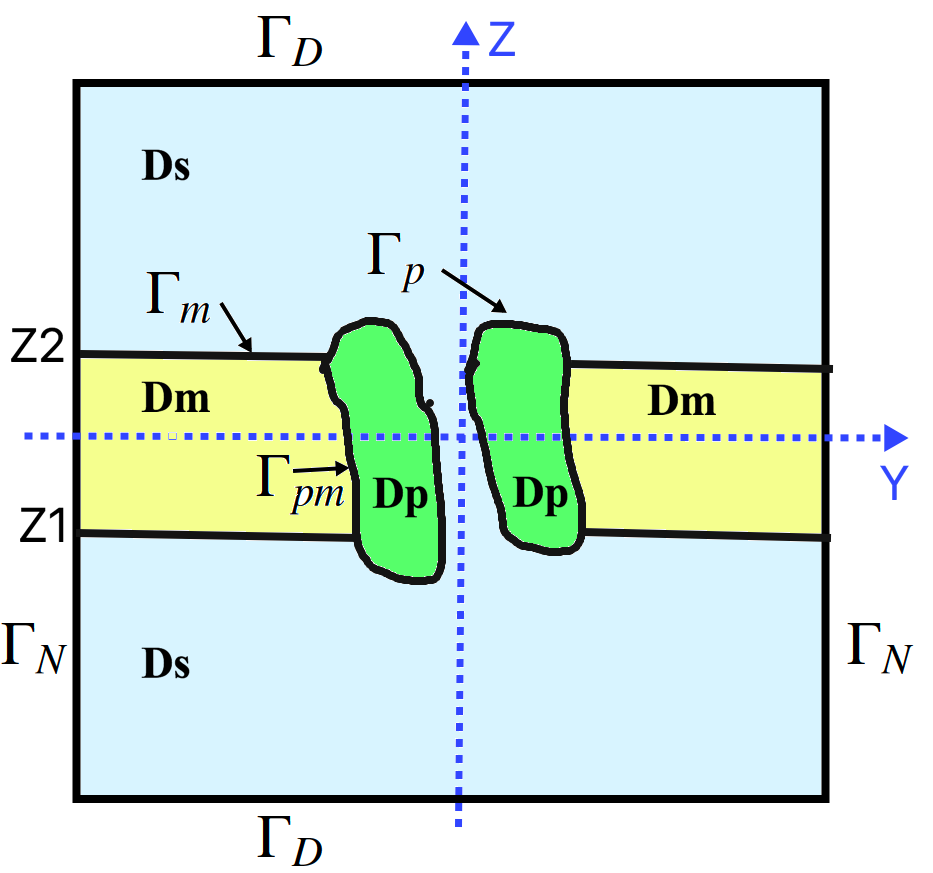}
  \end{center}
  \caption{\small An illustration of a box domain partition given in \eqref{DomainDecomp} and \eqref{DomainDecomp2}.}
  \label{boxdomain}
\end{figure} 

\section{A nonuniform size modified  Poisson-Boltzmann  ion channel  model}

We define a simulation box domain,   $\Omega$,  by
\begin{equation}\label{box-domain}
 \Omega =\{ (x,y,z) \; | \: L_{x1} < x < L_{x2}, \; L_{y1} < y < L_{y2}, \; L_{z1} < z < L_{z2} \;\},
\end{equation}
and split it into three open subdomains, $D_p$, $D_m$, and $D_s$ as follows: 
\begin{equation}\label{DomainDecomp}
\Omega = D_p\cup D_{m} \cup  D_s \cup \Gamma_p \cup \Gamma_m \cup \Gamma_{pm},
\end{equation}
where $L_{x1}, L_{x2}, L_{y1}, L_{y2}, L_{z1},$ and $L_{z2}$ are real numbers; 
$D_p$ is a protein region containing an ion channel protein molecule with $n_p$ atoms; $D_m$ is a  membrane region;  $D_s$  is a solvent region containing $n$ ionic species;  $\Gamma_p$ denotes an interface between $D_p$ and $D_s$; $\Gamma_m $ an interface between $D_{m}$ and $D_s$; and $\Gamma_{pm}$ an interface between $D_p$ and $D_m$.  We also split the boundary $\partial \Omega$ of $\Omega$  by
\begin{equation}\label{DomainDecomp2}
  \partial \Omega = \Gamma_D \cup \Gamma_N,
\end{equation}
where $\Gamma_D$ consists of the bottom and top surfaces of  $\Omega$  and $\Gamma_N$ consists of  the four side surfaces of $\Omega$. We further set the normal direction of the membrane surface in the $z$-axis direction so that  the membrane location can be determined by two numbers $Z1$ and $Z2$. An illustration of this setting and partitions \eqref{DomainDecomp} and \eqref{DomainDecomp2} is given in  Figure~\ref{boxdomain}. 

Let $c_i$ denote a concentration function of the $i$-th species in moles per liter (mol/L) and $u$ be a dimensionless electrostatic potential function of an electric field induced by atomic charges, ionic charges, and membrane charges. 
When a three-dimensional molecular structure of an ion channel protein and a mixture solution of $n$ ionic species  are  given, an atomic charge density function, $\rho_p$, and  an ionic charge density function,  $\rho_s$, can be estimated by
\[ \rho_p = e_c \sum\limits_{j=1}^{n_{p}}z_{j}  \delta_{\rr_{j}} \quad \mbox{ in } D_p, \qquad
   \rho_s(\rr) =  e_c \sum\limits_{i=1}^n Z_i c_i(\rr) \quad \mbox{ in } D_s,\]
where $e_{c}$ is  the elementary charge, $z_j$ and $\rr_j$ denote the atomic charge number and position vector of atom $j$, respectively, $Z_i$ is the charge number of ionic species $i$,  and $\delta_{\rr_{j}}$ is  the Dirac delta distribution at  $\rr_{j}$. In addition, a membrane surface charge density $\sigma$ is given in  $\mu$C/cm$^2$ to account for membrane charge effects. Based on the implicit solvent approach, the three regions $D_p$,  $D_{m}$, and $D_s$ are treated as  dielectric media with permittivity constants $\ep$, $\emm$, and $\es$, respectively. We then construct a Poisson dielectric boundary value problem for estimating  $u$ as follows:
\begin{equation}
  \label{Poisson-u}
\left\{ \begin{array}{rl}
  - \epsilon_p\Delta u( \rr)  =\alpha \sum\limits_{j=1}^{n_{p}}z_{j}  \delta_{\rr_{j}}, &        \rr \in D_p,  \\
   - \epsilon_m\Delta u( \rr)  =0, &         \rr \in D_{m},  \\
 - \epsilon_s \Delta u( \rr) =   \beta \sum\limits_{i=1}^n Z_i c_i(\rr), & \rr \in D_s,  \\  
u(\s^-) = u(\s^+), \quad \ep \frac{\partial u(\s^-)}{\partial \nn_p(\s)} = \es  \frac{\partial u(\s^+)}{\partial \nn_p(\s)},   & \s\in\Gamma_p,\\
  u(\s^-) = u(\s^+), \quad  \emm  \frac{\partial u(\s^-)}{\partial \nn_m(\s)} = \es \frac{\partial u(\s^+)}{\partial \nn_m(\s)} + \tau \sigma,  & \s\in\Gamma_m,\\
  u(\s^-) = u(\s^+), \quad  \ep  \frac{\partial u(\s^-)}{\partial \nn_p(\s)} = \emm \frac{\partial u(\s^+)}{\partial \nn_p(\s)},  & \s\in\Gamma_{pm},\\
   u(\s) = g(\s), & \s \in \Gamma_D,\\
\frac{\partial u( \s)}{\partial \nn_b(s)} = 0,  & s\in \Gamma_N,
\end{array}\right.
\end{equation}  
where $\nn_s$, $\nn_p$, $\nn_m$, and $\nn_b$ denote the unit outward normal directions of $D_s$, $D_{p}$,  $D_m$, and  $\Omega$,  respectively;  $\frac{\partial u(\s)}{\partial \nn(\s)}
$ denotes the directional derivative of $u$ along a unit  outside normal direction $\nn$ (say, $\nn=\nn_p$);   $u(\s^{\pm}) = \lim_{t\rightarrow 0^+} u(\s\pm t\nn(\s))$, which are the two sided limits along a direction $\nn$ of a region (say a protein region, $D_p$) from the inside and outside the region; $\alpha$,  $\beta$ and $\tau$ are three physical constants; and  $g$ is a boundary value function. Here 
the SI unit system has been used and the constants $\alpha$,  $\beta$ and $\tau$ are defined by
 \begin{equation}
\label{alpha-beta}
 \alpha = \frac{10^{10}e_{c}^{2}}{\ez k_{B}T}, \quad \beta = \frac{N_A e_{c}^{2}}{10^{17}\ez k_{B}T}, \quad \tau =  \frac{ 10^{-12} e_c}{\ez k_B T},
\end{equation} 
where $\ez$ is the permittivity of the vacuum, $k_B$ is the Boltzmann constant,  $T$ is the absolute temperature,  and $N_A$ is the  Avogadro number, which estimates  the number of ions per mole.  

From the Poisson problem \eqref{Poisson-u} it can be seen that different selections of ionic concentrations $c_i$ may lead to different electrostatic potentials $u$. To find an optimal selection, we define an electrostatic free energy functional, $F(c; u)$, by
\begin{equation}
\label{F-def}
   F(c; u) =F_{es}(c; u)+F_{id}(c)+F_{ex}(c) \qquad \mbox{with $\;\; c=(c_1,c_2,\ldots,c_n)$},
\end{equation}
where   $F_{es}$, $F_{id}$, and $F_{ex}$ denote the electrostatic, ideal gas, and excess energies, respectively, in the expressions 
\begin{eqnarray*}
F_{es}(c; u) &=&  \frac{k_B T}{2}\gamma  \sum\limits_{i=1}^n Z_i \int_{D_s } u c_{i} d\rr, \quad
F_{id}(c) =  k_{B}T\gamma  \sum\limits_{i=1}^n \int_{D_s } c_{i}\left( \ln \frac{c_{i}}{c_i^b} -1\right) d\mathbf{r},\\
 F_{ex}(c) &=& \frac{k_{B}T}{v_0}\int_{D_s } \left[1-\gamma  \sum\limits_{i=1}^n v_i c_i(\rr) \right] 
 \left[ \ln \left(1-\gamma  \sum\limits_{i=1}^n  v_i c_i(\rr)  \right) -1 \right] d\mathbf{r}. 
\end{eqnarray*}
Here, $v_i$ a volume of each ion of species $i$ in cubic angstroms  (\AA$^3$), $\gamma = 10^{-27} N_A$, $c_i^b$ is a bulk concentration,  and $v_0$ a size scaling parameter (e.g., $v_0=\min_iv_i$).  
Note that the sum $\gamma  \sum\limits_{i=1}^n v_i c_i(\rr)$ gives the ionic volume portion out of the solvent region. Hence, the difference $1 - \gamma  \sum\limits_{i=1}^n v_i c_i(\rr)$ is the portion of water solution volume, which should be positive, ensuring the definition of excess energy $F_{ex}$. These energy terms have been measured in energy units Joules. 

The constant $\gamma$ is a unit converter from mol/L to \AA$^3$.  In fact, under the SI unit system, a concentration is measured in the number of ions per \AA$^3$. But in practice, it is usually measured in mol/L as done in this paper. Hence, $\gamma$ is needed to convert from mol/L to  \AA$^3$ as shown below: 
\[  1 \mbox{ mol /L} =  10^{3} N_{A} /m^{3}= 10^{-27} N_{A}  /\mbox{\r{A}}^{3} = \gamma  /\mbox{\r{A}}^{3}, \]
where we have used the unit converters: 1 mol = 1 N$_A$, 1 L = $m^3/10^3$, and 1 m = $10^{10}$ \AA.
For  N$_A=6.02214129\times10^{23}$, $\gamma$ can be estimated as $\gamma \approx 6.022\times10^{-4}$. 

The free energy functional  $F$ of \eqref{F-def}  is a significant improvement of the conventional free energies  (e.g.,  \cite[eq. (1.4) or (3.4)]{BoLi2009a}, and  \cite[eq. (4)]{PhysRevLett.106.046102} or \cite[eq. (36)]{xie_PF2016}) since it discards the thermal de Broglie wavelengths, chemical potentials, a concentration of water molecules, and  terms $\ln(v_i c_i)$. 

The first  Fr{\'e}chet derivative $F^{\prime}$ of $F$ can be found in the expression 
\begin{equation*}
\label{fa-derivative}
 F^{\prime}(c; u)w= k_{B}T   \gamma  \sum\limits_{i=1}^{n}\int_{D_{s}}\Big[Z_i u+\ln\left(\frac{c_{i}}{c_{i}^{b}}\right) 
     -\frac{ v_{i}}{v_{0}} \ln \big(1-  \gamma  \sum\limits_{j=1}^n v_j c_j \big)  \Big] w_i(\rr)d\rr,
\end{equation*}
where $F^{\prime}(c; u)$ is a linear operator and $w$ denotes a test function vector with $w=(w_1,w_2,\ldots,w_n)$. From the above expression we can obtain the variation $\frac{\partial F(c; u)}{\partial c_i}$ of $F$ with respect to $c_i$ as follows: 
\begin{equation*}
  \label{mui_eq}
\frac{\partial F(c; u)}{\partial c_i} =  k_{B}T   \gamma  \left[ Z_i u + \ln\left(\frac{c_{i}}{c_{i}^{b}}\right)   -\frac{ v_{i}}{v_{0}} \ln \big(1-   \gamma  \sum\limits_{j=1}^n v_j c_j \big) \right], \quad i=1,2,\ldots, n.
\end{equation*}
Setting $\frac{\partial F(c; u)}{\partial c_i} =0$ gives the equation
$Z_i u + \ln\left(\frac{c_{i}}{c_{i}^{b}}\right)   -\frac{ v_{i}}{v_{0}} \ln \big(1-   \gamma  \sum\limits_{j=1}^n v_j c_j \big) =0$, which can be reformulated as 
\begin{equation}
  \label{sizeModifiedBoltzmann}
  c_{i}(\rr) - {c}_{i}^b  \left[1-  \gamma  \sum\limits_{j=1}^n  v_j c_j(\rr)  \right]^{\frac{ v_i}{v_0}} e^{-Z_{i} u(\rr)} =0,
  \quad  i=1,2,\ldots,n.
\end{equation}  
Obviously, the above nonlinear algebraic equations give the necessary conditions that an optimal $c$ must satisfy. They also can be regarded as the ion size constraint conditions during a search process for an optimal $u$. Because  the nonlinear equations  \eqref{sizeModifiedBoltzmann} and the Poisson problem \eqref{Poisson-u} are dependent each other, they must be combined together as a nonlinear system in terms of $c$ and $u$.  A solution of this system gives an optimal $c$ and  an optimal $u$ in the sense of minimizing the electrostatic free energy \eqref{F-def}. This nonlinear system defines a nonuniform size modified Poisson-Boltzmann ion channel (nuSMPBIC) mode. For clarity, we restate it as follows:
\begin{equation}
  \label{nuSMPBic}
\left\{ \begin{array}{rl}
 c_{i}(\rr) - {c}_{i}^b  \left[1-  \gamma  \sum\limits_{j=1}^n  v_j c_j(\rr)  \right]^{\frac{ v_i}{v_0}} e^{-Z_{i} u(\rr)} =0,
  &  \rr\in D_{s}, \quad  i=1,2,\ldots,n,\\
  - \epsilon_p\Delta u( \rr)  =\alpha \sum\limits_{j=1}^{n_{p}}z_{j}  \delta_{\rr_{j}}, &        \rr \in D_p,  \\
   - \epsilon_m\Delta u( \rr)  =0, &         \rr \in D_{m},  \\
  \epsilon_s \Delta u( \rr) -   \beta \sum\limits_{i=1}^n Z_i c_i(\rr) = 0, & \rr \in D_s,  \\  
u(\s^-) = u(\s^+), \quad \ep \frac{\partial u(\s^-)}{\partial \nn_p(\s)} = \es  \frac{\partial u(\s^+)}{\partial \nn_p(\s)},   & \s\in\Gamma_p,\\
  u(\s^-) = u(\s^+), \quad  \emm  \frac{\partial u(\s^-)}{\partial \nn_m(\s)} = \es \frac{\partial u(\s^+)}{\partial \nn_m(\s)} + \tau \sigma,  & \s\in\Gamma_m,\\
  u(\s^-) = u(\s^+), \quad  \ep  \frac{\partial u(\s^-)}{\partial \nn_p(\s)} = \emm \frac{\partial u(\s^+)}{\partial \nn_p(\s)},  & \s\in\Gamma_{pm},\\
   u(\s) = g(\s), & \s \in \Gamma_D,\\
\frac{\partial u( \s)}{\partial \nn_b(s)} = 0,  & s\in \Gamma_N.
\end{array}\right.
\end{equation}  

In physics, the Neumann boundary condition on the four  side surface $\Gamma_{N}$ reflects the fact that none of the charges enter the box domain $\Omega$ from $\Gamma_{N}$. To mimic a voltage across the membrane, we set the boundary function $g$ as a piecewise function: 
\begin{equation}
    \label{g_def}
 g(\s) = \left\{\begin{array}{rl} u_b & \mbox{at  $z= L_{z1}$ (bottom surface)}, \\
                                               u_t  & \mbox{at  $z=L_{z2}$ (top surface)},
 \end{array}\right.                                      
\end{equation}
where $u_{b}$ and $u_{t}$ denote two electrostatic potential values. 

The nuSMPBIC model  involves Dirac delta distributions, two physical domains (the solvent domain $D_s$ and the box domain $\Omega$), complicated interface conditions, mixed boundary value conditions (i.e., a Dirichlet boundary condition on $\Gamma_D$ and a Neumann boundary condition on $\Gamma_N$), and the membrane surface charge density $\sigma$. The  Dirac delta distributions  cause the potential function $u$ strongly singular  while the nonlinear algebraic equations cause the ionic concentrations $c_i$ exponentially nonlinear. Hence, the nuSMPBIC model \eqref{nuSMPBic} is very difficult to solve numerically. New numerical techniques are needed to develop an effective nuSMPBIC finite element solver.  

\section{A submodel partition of  the nuSMPBIC model}
One major difficulty in the numerical solution of the nuSMPBIC model comes from the solution singularity caused by the  Dirac-delta distributions $\delta_{\rr_{j}}$. Following what are done in \cite{xiePBE2013,SMPBEic2019}, in this section, we partition the nuSMPBIC model into three submodels, called Models 1, 2, and 3, to overcome the singularity difficulty.

Model 1 is defined by the Poisson equation over the whole space $\R^3$,
\begin{equation}
- \ep \Delta G(\rr) =\alpha \sum\limits_{j=1}^{n_{p}}z_{j}  \delta_{\rr_{j}},  \quad  \rr \in \R^3, \quad G(\rr) \rightarrow 0 \mbox{ as } |\rr| \rightarrow \infty,
\label{Poisson-Dp2}
\end{equation} 
whose solution $G$ gives an electrostatic potential induced by atomic charges from an ion channel protein. The analytical solution $G$ of Model~1 and its gradient vector $\nabla G(\s)$ can be found in the algebraic expressions \cite{xiePBE2013}
\begin{equation}
\label{G-def}
G(\rr)= \frac{\alpha}{4\pi\ep} \sum_{j=1}^{n_p}\frac{z_{j}}{| \rr-\rr_j |},  \qquad  \nabla G(\rr) = -\frac{\alpha}{4\pi\ep}  \sum_{j=1}^{n_p}z_{j}\frac{(\rr-\rr_j)}{|\rr-\rr_j|^3}     \qquad     \forall \rr \neq \rr_j.
\end{equation}  

Model 2 is defined by the linear interface boundary value problem:  
 \begin{equation}\label{Psi-channel}
\left\{\begin{array}{lll}
 \Delta \Psi(\rr)=0, \qquad \rr\in D_{m}\cup D_p\cup D_s, & \\
\Psi(\s^-)=\Psi(\s^+),  \quad
 \ep \frac{\partial \Psi(\s^-)}{\partial\nn_p(\s)}=\es \frac{\partial \Psi(\s^+)}{\partial\nn_p(\s)}+(\es-\ep)\frac{\partial G(\s)}{\partial \nn_p(\s)}, &  \s\in\Gamma_p, \\
 \Psi(\s^-)=\Psi(\s^+),  \quad
 \emm \frac{\partial \Psi(\s^-)}{\partial\nn_m(\s)}=\es \frac{\partial \Psi(\s^+)}{\partial\nn_m(\s)}
 +(\es-\emm)\frac{\partial G(\s)}{\partial \nn_m(\s)} + \tau \sigma, &  \s\in\Gamma_m, \\
  \Psi(\s^-)=\Psi(\s^+),  \quad
 \ep \frac{\partial \Psi(\s^-)}{\partial\nn_p(\s)}=\emm \frac{\partial \Psi(\s^+)}{\partial\nn_p(\s)}
 +(\emm-\ep)\frac{\partial G(\s)}{\partial \nn_p(\s)}, &  \s\in\Gamma_{pm}, \\
 \Psi(\s) ={g}(\s) -G(\s),   &  \s\in \Gamma_D,  \\
  \frac{\partial \Psi(\s)}{\partial\nn_b(\s)} = - \frac{\partial G(\s)}{\partial \nn_b(\s)}, & \s \in \Gamma_{N},
\end{array}\right.
\end{equation}
Clearly, both Models 1 and 2 are independent of ionic concentration $c_i$. Hence, they can be solved prior to a search for ionic concentrations $c_i$. 

With Models 1 and 2, we can simplify the nuSMPBIC model \eqref{nuSMPBic} into Model 3 --- a nonlinear system of an electrostatic potential, $\Phit$, (induced purely by ionic charges) and ionic concentration functions $c_i$  as follows:
\begin{equation}
\label{smpbeic-def}
\left\{\begin{array}{ll}
 c_{i}(\rr) - {c}_{i}^b  \left[1-  \gamma  \sum\limits_{j=1}^n  v_j c_j(\rr)  \right]^{\frac{ v_i}{v_0}} e^{-Z_{i} [G(\rr) + \Psi(\rr) + \Phit(\rr)]} =0, & \rr\in D_s,\; i=1,2,\ldots, n,\\
\Delta \Phit(\rr)=0, &\rr\in D_{m}\cup D_p,\\
\es\Delta  \Phit(\rr) + \beta \sum\limits_{i=1}^n Z_{i}  c_i(\rr) = 0, & \rr\in D_s,\\
  \Phit(\s^+)= \Phit(\s^-), \quad  \es \frac{\partial \Phit(\s^+)}{\partial\nn_p(\s)}=\ep \frac{\partial \Phit(\s^-)}{\partial\nn_p(\s)}, &  \s\in\Gamma_p,\\
   \Phit(\s^+)= \Phit(\s^-), \quad  \es \frac{\partial \Phit(\s^+)}{\partial\nn_m(\s)}=\emm \frac{\partial \Phit(\s^-)}{\partial\nn_m(\s)}, &  \s\in\Gamma_m,\\
    \Phit(\s^-)=\Phit(\s^+),  \quad
 \ep \frac{\partial \Phit(\s^-)}{\partial\nn_p(\s)}=\emm \frac{\partial \Phit(\s^+)}{\partial\nn_p(\s)}, &  \s\in\Gamma_{pm}, \\
 \Phit(\s) = 0,  & \s \in \Gamma_D,  \\
 \frac{\partial \Phit(\s)}{\partial \nn_b(s)} = 0,  & s\in \Gamma_N.
 \end{array}\right.
\end{equation}

After solving Models 2 and 3, we construct the electrostatic potential function $u$ by the formula
\begin{equation}\label{solutionSplit}
u(\rr)=G(\rr) + \Psi(\rr) + \Phit(\rr) \quad \quad   \forall \rr\in \Omega.  
\end{equation}

Since  $G$ collects all the singularity points of $u$, both Models 2 and 3 can be much easier to solve numerically than the original nuSMPBIC model. Consequently, the complexity of solving the nuSMPBIC model is sharply reduced. 
An efficient finite element method for solving Model 2 has been reported in \cite{SMPBEic2019,Xie4PNPicNeumann2020}. Hence, we only need to develop a finite element method for solving Model 3 as done in the next section.

\section{A damped two-block iterative method for solving Model 3}
We start with a construction of two finite element meshes --- an interface fitted irregular tetrahedral mesh, $\Omega_h$, of a box domain $\Omega$ and  a tetrahedral mesh, $D_{s,h}$, of  a solvent domain $D_s$. We then use them construct  two linear Lagrange finite element function spaces, denoted by $\mathcal{U}$ and $\mathcal{V}$, as two finite dimensional subspaces of the  Sobolev function spaces $H^1(\Omega)$ and $H^1(D_s)$ \cite{adams2003sobolev}, respectively. We further construct a restriction operator, $ {\cal R}: \mathcal{U} \rightarrow \mathcal{V}$, and a prolongation operator, ${\cal P}: \mathcal{V} \rightarrow \mathcal{U}$, such that ${\cal R} u \in \mathcal{V}$ for any $u\in \mathcal{U}$ and ${\cal P} c_i \in \mathcal{U}$ for any $c_i\in \mathcal{V}$.
With these spaces and operators, we can obtain a finite element approximation of Model 3 as follows: {Find  $\Phit \in \mathcal{U}_0$ and $c_i\in \mathcal{V}$ for $i=1,2,\ldots,n$ such that}
\begin{subequations}
\label{transformed_system}
\begin{eqnarray}
   \label{EQ4C}
    c_{i}(\rr) - {c}_{i}^b   \left[1-  \gamma  \sum\limits_{j=1}^n  v_j c_j(\rr)  \right]^{\frac{ v_i}{v_0}} e^{-Z_{i} {\cal R} [G(\rr) + \Psi(\rr) + \Phit(\rr)]} =0, 
       &  \rr\in D_{s,h}, \quad  i=1,2,\ldots,n,\\
      \label{EQ4Phit}
        a(\Phit,v)   - \beta \sum\limits_{j=1}^{n}Z_{j} \int_{D_{s}} {\cal P} {c}_{j}(\rr) v d\rr   = 0 &  \forall v  \in \mathcal{U}_0,
\end{eqnarray}
\end{subequations}
where $\mathcal{U}_0 =  \{ u \in \mathcal{U} \; | \:  u  =0 \mbox{ on } \Gamma_D \}$, which is a subspace of $\mathcal{U}$, and  $a(\cdot, \cdot)$  is a bilinear form defined by
 \begin{equation}
\label{WeakForm1} 
 a(u, v) =  \ep\int_{D_{p}}\nabla u \cdot \nabla v d\rr +  \emm\int_{D_{m}}\nabla u \cdot \nabla v d\rr + \es \int_{D_{s}} \nabla u \cdot \nabla v d\rr, \quad u,v\in \mathcal{U}.
\end{equation}

We now present the damped two-block iterative method for solving the finite element system \eqref{transformed_system}. Here 
we divide the unknown functions of \eqref{transformed_system} into two blocks with Block 1 containing all the ionic concentration functions $c_i$ and Block 2 containing $\Phit$ only. With vector $c=(c_1,c_2, \ldots, c_n)$,  we rewrite \eqref{EQ4C} in the vector equation
\begin{equation}
\label{System4c}
      F(c(\rr), \Phit(\rr)) = \mathbf{0}  \quad \forall \rr  \in D_{s,h}, 
\end{equation}
where $F =\left(f_{1}, f_{2}, \ldots, f_n \right)$ with the $i$-th component function $f_i$ being defined by
\begin{equation}
 \label{System4c2}   
f_i(c(\rr), \Phit(\rr)) =  c_{i}(\rr) - {c}_{i}^b \left[1-  \gamma  \sum\limits_{j=1}^n  v_j c_j(\rr)  \right]^{\frac{ v_i}{v_0}} e^{-Z_{i} {\cal R} [G(\rr) + \Psi(\rr) + \Phit(\rr)]},\quad \rr \in D_{s,h}, \quad  i=1,2,\ldots,n.
\end{equation}
Thus, \eqref{transformed_system} has been rewritten in a two block form --- a system of \eqref{System4c} and \eqref{EQ4Phit}.

Let $(c^k, \Phit^k)$ denote the $k$-th iterate of the damped two-block iterative method with $c^k=(c_1^k,c_2^k,\ldots, c_n^k)$. When an initial iterate, $(c^0, \Phit^0)$, is given, we define the damped two-block iterative method  as follows:
\begin{subequations}
\label{two-block-outloop}
\begin{eqnarray}
\label{cj-iterate}
    {c}^{k+1}(\rr) &=&  {c}^{k}(\rr) + \omega \left[ {p}(\rr)  - {c}^{k}(\rr) \right], \quad \rr \in D_{s,h},\\
    \label{phit-iterate}
     \Phit^{k+1}(\rr) &=& \Phit^{k}(\rr) +  \omega \left[{q}(\rr) - \Phit^{k}(\rr) \right], \quad \rr \in \Omega,  \quad k = 0, 1, 2, \ldots,
\end{eqnarray}
\end{subequations}
where $\omega$ is a damping parameter between 0 and 1,  $p=(p_1,p_2,\ldots, p_n)$ is a solution of the nonlinear algebraic system 
\begin{equation}
 \label{System4p}
        F(p(\rr), \Phit^k(\rr)) = \mathbf{0}, \quad \rr \in D_{s,h},
\end{equation}
and $q$ is a solution  of the linear finite element variational problem:  \mbox{Find $q \in  \mathcal{U}_0$ such that}    
\begin{equation}
 \label{eq4q}
   a(q, v)   = \beta \sum\limits_{j=1}^{n}Z_{j} \int_{D_{s}} {\cal P} {c}_{j}^{k+1}(\rr) v d\rr  \quad  \forall v  \in \mathcal{U}_0.
\end{equation}

In the implementation, we use the following iteration termination rules:
\begin{equation}
    \label{Ite-stop}
    \|   \Phit^{k+1}  -  \Phit^{k}  \|_{\Omega} < \epsilon,  \quad  \max_{1\leq i\leq n} \|   {c}^{k+1}_i -  {c}^{k}_i  \|_{D_s} < \epsilon,
      \quad \mbox{and }  \quad    R({c}^{k+1}, \Phit^{k+1}) < \epsilon,
\end{equation} 
where $\epsilon$ is a tolerance (by default, $\epsilon=10^{-4}$),  $ \| \cdot \|_{\Omega}$ and  $ \| \cdot  \|_{D_s}$ denote the norms of function spaces $L_2(\Omega)$ and $L_2(D_s)$, respectively, and $R({c}, \Phit)$ denotes a residual error of the nonlinear algebraic system \eqref{System4c} as below:
\begin{equation}
    \label{constrianResidual}
       R({c}, \Phit) = \frac{1}{N_h} \max_{1\leq i\leq n} \left( \sum_{\mu=1}^{N_h} \big| f_i(c(\rr^{(\mu)}), \Phit(\rr^{(\mu)}) \big |^2 \right)^{1/2}.
\end{equation} 
Here $\rr^{(\mu)}$ denotes the $\mu$-th mesh point of $D_{s,h}$ and $N_h$ is the total number of mesh points. Clearly, 
 $R({c}, \Phit)=0$ if and only if $F(c(\rr), \Phit(\rr)) = \mathbf{0}$ for $\rr \in D_{s}$. 
We stop the iterative process \eqref{two-block-outloop} and output the $(k+1)$-th iterate $( {c}^{k+1}, \Phit^{k+1})$ as an approximate solution  $(c,\Phit)$ of Model 3 whenever  the iteration termination rules \eqref{Ite-stop} are satisfied. 

\section{A Newton iterative method for solving a nonlinear algebraic system of Block 1 }
In this section, we present a novel Newton iterative scheme for solving the nonlinear algebraic system \eqref{System4p} of Block 1. It is this efficient scheme that turns the damped two-block iterative method \eqref{two-block-outloop} into an efficient algorithm for solving Model 3.

A construction of such a Newton iterative scheme is motivated from one basic property of  a linear finite element function. That is,  a linear finite element function is determined uniquely by its mesh node values. In other words, to search for a unknown linear finite element function, we only need to determine its mesh node values. According to this property, we set $\rr= \rr^{(\mu)}$ in the equation of \eqref{System4p} to produce $N_h$ small nonlinear systems as follows:
\begin{equation}
 \label{System4cPhit}
   \bar{F}(\xi_\mu) = \mathbf{0}, \quad  \mu =1,2,\ldots, N_h,
\end{equation}
where $\xi_\mu = ( \xi_{1,\mu}, \xi_{2,\mu}, \ldots, \xi_{n,\mu})$ with $\xi_{i,\mu}$ denoting the mesh node value $p_i(\rr^{(\mu)})$ of $p_i$  and $\bar{F} =\left(\bar{f}_{1}, \bar{f}_{2}, \ldots, \bar{f}_n \right)$ with $\bar{f}_i$ being defined by
\[  \bar{f}_i(\xi_\mu) =  \xi_{i,\mu} - {c}_{i}^b \left[1-  \gamma  \sum\limits_{j=1}^n  v_j \xi_{j,\mu}  \right]^{\frac{ v_i}{v_0}} 
     e^{-Z_{i} u_{k,\mu}}   \qquad \mbox{with }\quad u_{k,\mu} = G(\rr^{(\mu)}) + \Psi(\rr^{(\mu)}) + \Phit^k(\rr^{(\mu)}).\]
Here $u_{k,\mu}$ is a known value but may be too large to cause a numerical overflow problem during an iterative process. To avoid the overflow problem, for a given allowable upper bound $M$, we modify $u_{k,\mu}$ as $- M/Z_i$ whenever
$  -Z_i u_{k,\mu} \geq M$. By default, we set $M=45$. 

Clearly, the $N_h$ nonlinear systems of \eqref{System4cPhit} are independent each other. Hence, they can be solved  one-by-one independently to produce $nN_h$ mesh node values $p_i(\rr^{(\mu)})$ for $i=1,2,\ldots,n$ and $\mu=1,2,\ldots, N_h$. We then use them to derive a solution of the nonlinear algebraic system~\eqref{System4p}.

We now construct a Newton iterative scheme for solving each small nonlinear system of \eqref{System4cPhit}. 
Let $\xi_{\mu}^j$ denote the $j$-th iterate of the  Newton iterative scheme. When an initial iterate $\xi_{\mu}^0$ is given (by default, $\xi_{\mu}^0= c^{k}$), we define the Newton iterative scheme by 
\begin{equation}
    \label{Newton4p}
   \xi_{\mu}^{j+1} = \xi_{\mu}^j  + \Upsilon_j,   \quad j=0, 1, 2, \ldots,
\end{equation} 
where $\Upsilon_j$ is a solution of the Newton equation  
\begin{equation}
    \label{Newton4r}
  J(  \xi_{\mu}^j)  \Upsilon_j = - \bar{F}(\xi_{\mu}^j). 
\end{equation} 
Here $J$ denotes a $n\times n$ Jacobian matrix of $\bar{F}$ with the $(i,j)$-th entry being the partial derivatives $\partial \bar{f}_{i}/\partial \xi_{j,\mu}$ for $i,j=1,2, \ldots,n$ as follows: 
\begin{equation}
\label{JacobianEntries}
\frac{\partial \bar{f}_{i}}{\partial \xi_{j,\mu}} = \left\{\begin{array}{ll}
 1 + \gamma \frac{v_i^2}{v_0} {c}_{i}^{b} e^{-Z_{i} u_{k,\mu}}   \left[1-  \gamma  \sum\limits_{j=1}^n  v_j \xi_{j,\mu}
       \right]^{\frac{ v_i}{v_0} - 1}, & j = i,\\
 \gamma \frac{v_i v_j}{v_0} {c}_{i}^{b} e^{-Z_{i} u_{k,\mu}}  \left[1-  \gamma  \sum\limits_{j=1}^n  v_j \xi_{j,\mu}
       \right]^{\frac{ v_i}{v_0} - 1}, & j \neq i.
\end{array}\right.
\end{equation}
Each Newton equation \eqref{Newton4r} can be solved directly by the Gaussian elimination method since  $n$ is small (say 2, 3, or 4).  When $\| \Upsilon_j \| < \epsilon $ (by default $\epsilon=10^{-8}$), we output $\xi_{\mu}^{ j+1}$ as a numerical solution of \eqref{System4cPhit}.

\section{A good initial iterate for the damped two-block\\ iterative method}
When all the ion sizes $v_i$ are set to be equal to the average size $\bar{v} = \frac{1}{n} \sum_{i=1}^n v_i $, we can solve the nonlinear algebraic system \eqref{sizeModifiedBoltzmann} for $c$ to derive an analytical expression of $c_i$ as follows:
\begin{equation}
 \label{c-selection}
  c_i(\rr) = \frac{c_i^b e^{-Z_i u(\rr)} }{1 +  \gamma  \frac{\bar{v}^2}{v_0}   \sum\limits_{j=1}^n c_j^b e^{-Z_j u(\rr) }}, \quad \rr \in D_s, \quad i=1,2,\ldots, n.
\end{equation}
Substituting the above expressions to the finite element equation \eqref{EQ4Phit} and using  the formula $u= G+\Psi + \Phit$, we can obtain the following nonlinear finite element variational problem:  {Find $\Phit \in\mathcal{U}_{0}$ such that  } 
\begin{equation}
 \label{uSMPhit-def}
         a(\Phit,v)   - \beta  \int_{D_{s}} \frac{ \sum\limits_{i=1}^{n}Z_{i} c_i^b e^{-Z_i (G+\Psi + \Phit )}}{1
+  \gamma  \frac{\bar{v}^2}{v_0}   \sum\limits_{i=1}^n c_i^b e^{-Z_i (G+\Psi + \Phit )}} v d\rr   = 0 
 \quad  \forall v \in \mathcal{U}_{0},
\end{equation}

Using the Taylor expansion $e^x = 1+x + O(x^2)$ and the electroneutrality condition $\sum_{i=1}^n Z_ic_i^b = 0$, we can linearize the nonlinear problem \eqref{uSMPhit-def} as the following linear finite element equation:  {Find $\Phit \in \mathcal{U}_{0}$ such that } 
\begin{equation}
    \label{Phit-iterate0}
 a(\Phit, v) +  \bar{\beta} \int_{D_{s}} \Phit v d\rr  = -  \bar{\beta} \int_{D_{s}} (G+\Psi) v d\rr    \quad  \forall v \in \mathcal{U}_{0}.
\end{equation}
where  $\bar{\beta}$ is defined by
\[ \bar{\beta} =  \frac{\beta \sum\limits_{i=1}^{n} Z_i^2 c_i^b}{1+\gamma \frac{\bar{v}^2}{v_0}  \sum\limits_{i=1}^{n} c_i^b}.\] 
A  solution of the above linear equation is selected as an  initial iterate $\Phit^{(0)}$ of the damped two-block iterative method  \eqref{two-block-outloop}.  We then use the formula \eqref{c-selection} to get the initial iterates $c_i^0$  by
\begin{equation}
 \label{c0-selection}
  c_i^0(\rr) = \frac{c_i^b e^{-Z_i [ G(\rr)+\Psi(\rr) + \Phit^0(\rr)]} }
                           {1 +  \gamma  \frac{\bar{v}^2}{v_0}   \sum\limits_{j=1}^n c_j^b e^{-Z_j [ G(\rr)+\Psi(\rr) + \Phit^0(\rr)] }}, 
                           \quad \rr \in D_s, \quad i=1,2,\ldots, n.
\end{equation}

 \section{Program package and numerical results}

We have presented a nuSMPBIC finite element solver in the above sections. For clarity, we summarize it in Algorithm 1. We then implemented it in Python and Fortran as a software package based on the state-of-the-art finite element library from the FEniCS project \cite{fenics-book} and the Poisson-Boltzmann finite element program packages reported in  \cite{xiePBE2013, nuSMPBE2017,SMPBEic2019}. Here the linear interface boundary value problem \eqref{Psi-channel} of Model 2 is solved by the efficient finite element method reported in \cite{SMPBEic2019,Xie4PNPicNeumann2020} while other related linear finite element equations are solved approximately by a  generalized minimal residual method using incomplete LU preconditioning (GMRES-ILU) with the absolute and relative residual error tolerances being  $10^{-6}$ by default.

\vspace{1mm}
{\bf Algorithm 1 (The nuSMPBIC finite element solver).} 
{\em Let  ${\cal U}$ and ${\cal V}$  be  linear Lagrange finite element spaces of $H^1(\Omega)$ and $H^1(Ds)$, respectively, and  $(c^k, \Phit^k)$ with $c^k=(c_1^k,c_2^k,\ldots, c_n^k)$ be the $k$-th  iterate of the damped two-block iterative method \eqref{two-block-outloop}.  A finite element solution $(u, c)$ of the nuSMPBIC model \eqref{nuSMPBic} with $c=(c_1,c_2,\ldots, c_n)$ is calculated in the following four steps:
\begin{enumerate}
\item[Step 1.]  Initialization:
   \begin{enumerate}
   \item  Calculate G and $\nabla G$ via \eqref{G-def} on ${\cal U}$.
   \item  Solve Model 2 for  $\Psi$ on ${\cal U}$.
   \item Calculate $\Phit^0$ by solving  the linear finite element equation \eqref{Phit-iterate0} on ${\cal U}$ and $c^0$ by formula \eqref{c0-selection}.
\item Set $k=0$.
\end{enumerate}
\item[Step 2.]  Calculate $(c^{k+1}, \Phit^{k+1})$  by the damped two-block iterative method \eqref{two-block-outloop}. Here the nonlinear algebraic system \eqref{System4p} is solved by the Newton iterative scheme \eqref{Newton4p}.
\item[Step 3.]  Check the convergence: If \eqref{Ite-stop} holds, set  $(c^{k+1}, \Phit^{k+1})$ as 
an approximate solution of Model 3 and go to Step 4; otherwise, increase $k$ by 1 and go back to Step~2.
\item[Step 4.]  Construct u by the solution decomposition \eqref{solutionSplit}: $u(\rr)=G(\rr) + \Psi(\rr) + \Phit(\rr) \quad   \forall \rr\in \Omega$.
\end{enumerate}
}

\begin{figure}[t]
        \centering
         \begin{subfigure}[b]{0.4\textwidth}
                \centering
                \includegraphics[width=\textwidth]{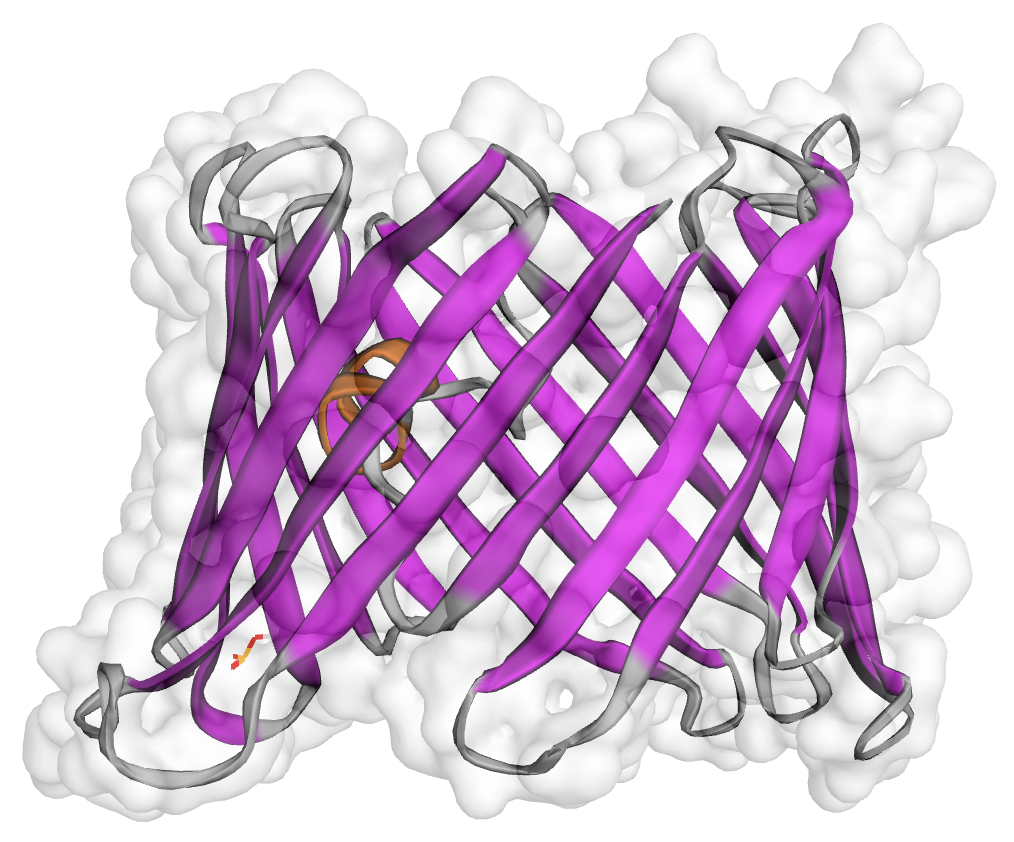}
                \vspace{1mm}
                \caption{Side view of mVDAC1}
        \end{subfigure}  
        \qquad \qquad \qquad
         \begin{subfigure}[b]{0.4\textwidth}
                \centering
                \includegraphics[width=\textwidth]{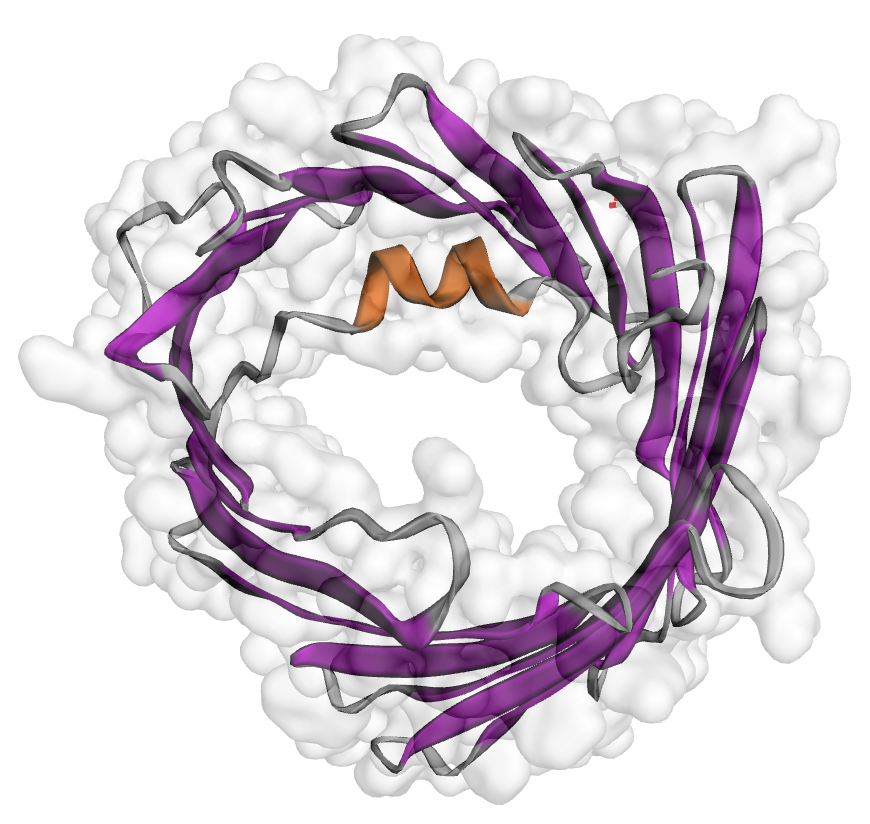}
                \caption{Top view of mVDAC1}
        \end{subfigure}  
              \caption{A crystallographic three-dimensional molecular structure of mVDAC1 (PDB ID: 3EMN) depicted in cartoon representations. Here the Van der Waals volume (a volume occupied by all the individual atomic balls of mVDAC1) is also displayed in grey color. 
              }        
\label{structure_3EMN}          

\end{figure}
\begin{figure}
        \centering
         \begin{subfigure}[b]{0.3\textwidth}
                \centering
                \includegraphics[width=\textwidth]{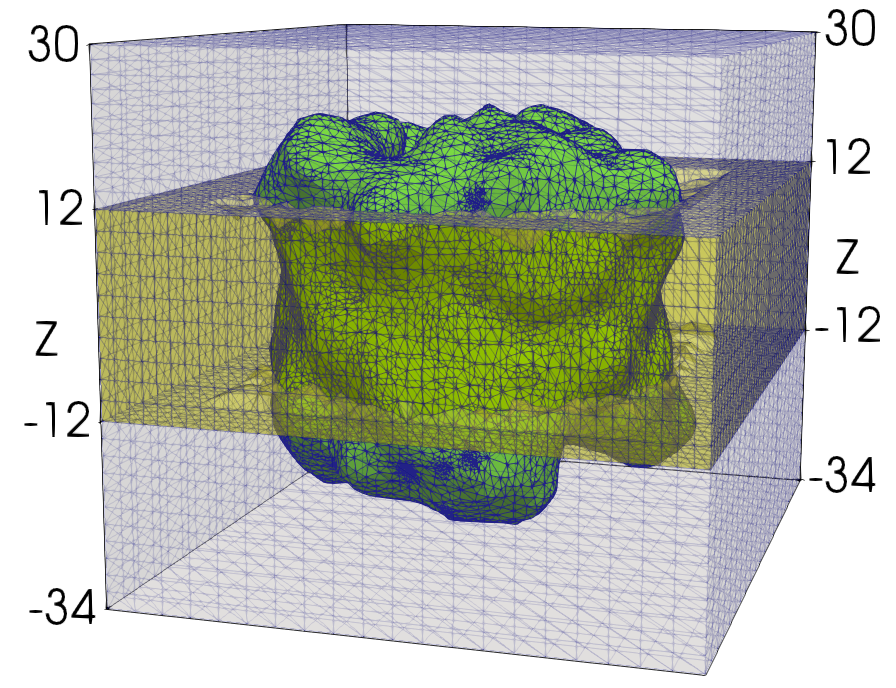}
                \caption{Box domain mesh $\Omega_h$}
        \end{subfigure}  
        \quad 
         \begin{subfigure}[b]{0.3\textwidth}
                \centering
                \includegraphics[width=\textwidth]{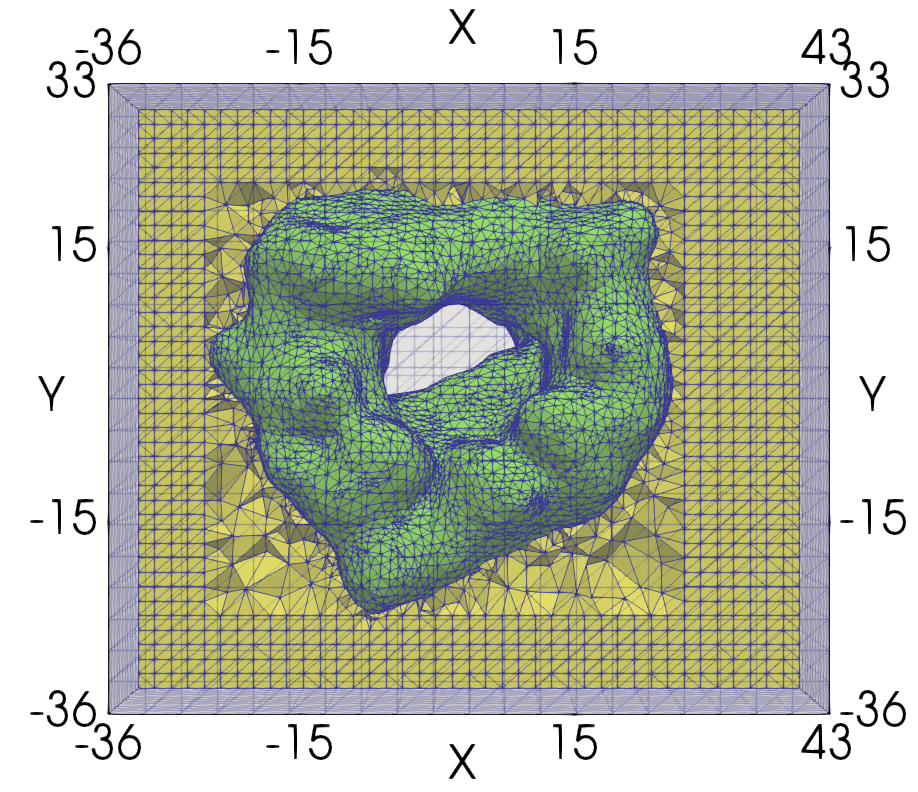}
                \caption{Top view of  $\Omega_h$}
        \end{subfigure}  
        \quad 
         \begin{subfigure}[b]{0.3\textwidth}
                \centering
                \includegraphics[width=\textwidth]{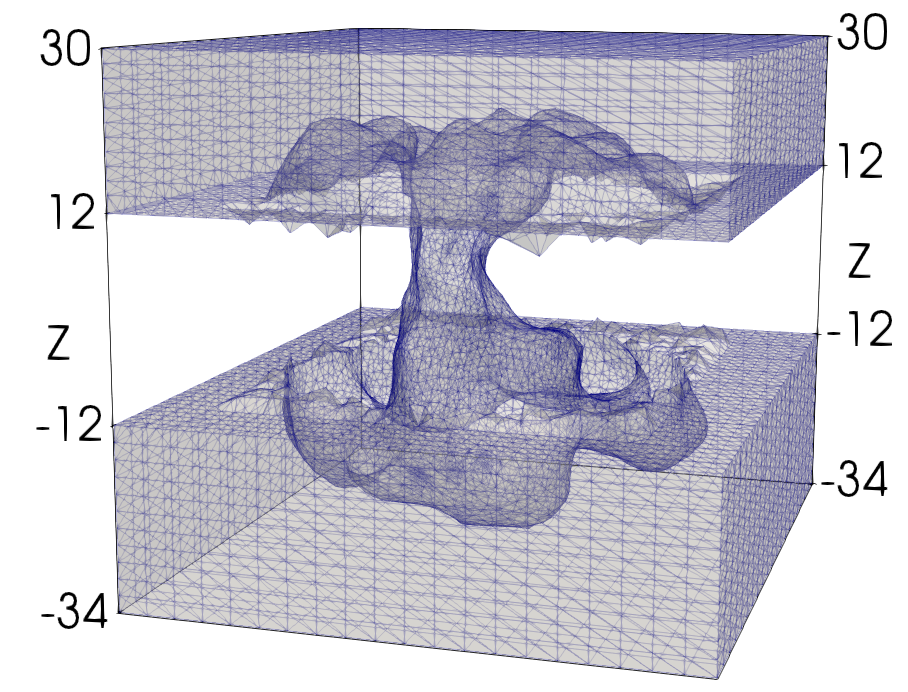}
                \caption{Solvent region mesh $D_{s,h}$ }
        \end{subfigure}
       \caption{(a, b) Two views of an  interface fitted irregular tetrahedral mesh $\Omega_h$ of  a box domain $\Omega$. (c) An  irregular tetrahedral mesh $D_{s,h}$ of  the solvent region $D_s$ extracted from $\Omega_h$. Here the meshes of the membrane region $D_m$ and protein region $D_p$ are colored in yellow and green and the two solvent compartments above and below the membrane belong to the cytoplasm and  intermembrane space, respectively.}        
\label{mesh_3EMN}          
\end{figure}

To demonstrate the convergence of the damped two-block iterative method for solving Model 3 and the performance of our nuSMPBIC finite element package, we did numerical tests on a murine voltage-dependent anion channel 1 (mVDAC1)  \cite{ujwal2008crystal} in a mixture of  0.1 mole KNO$_3$ (potassium nitrate) and 0.1 mole NaCl (table salt). Here the four ionic species Cl$^-$, NO$_3^-$, K$^+$,  and Na$^+$ were ordered from 1 to 4 for their concentration functions $c_i$, bulk concentrations $ c_i^b=0.1$, and charge numbers  $Z_1=-1, Z_2=-1, Z_3=1$, and $Z_4=1$. We treated each ion as a ball to estimate $v_i$ via the ball volume formula $v_i=4\pi r_i^3/3$ with $r_i$ denoting an ionic radius in \AA. From the website {\em https://bionumbers.hms.harvard.edu/bionumber.aspx?\&id=108517} we got 
\[  r_1=1.81, \; r_2=2.64,\; r_3=1.33,\; r_4=0.95. \]
Using them, we obtained the four ion sizes $v_i$ and their average volume $\hat{v} = (v_1+v_2+v_3+v_4)/4$ as follows:
\[ v_1= 24.8384, \quad v_2=77.0727, \quad  v_3=9.8547,  \quad v_4=3.5914, \quad \hat{v} = 28.8393. \]
This mixture is a good selection for us to demonstrate the importance of considering distinct ion sizes since it contains two anions with the same charge number $-1$ (Cl$^-$ and NO$_3^-$) and two cations with the same charge number $+1$ (K$^+$ and Na$^+$) and have four  significantly different  ion sizes.

As the main conduit on the outer mitochondrial membrane for the entry and exit of ions and metabolites between the cytosol and the mitochondria, the mVDAC1 has various ionic species with very different ion sizes within its channel pore, making it a good test case for us to validate the nuSMPBEic model not mention its anion-selectivity property and complex molecular structure \cite{ujwal2008crystal}. 
We used a crystallographic three-dimensional molecular structure of  mVDAC1 from  the Orientations of Proteins in Membranes (OPM) database ( {\em https://opm.phar.umich.edu}), instead of the Protein Data Bank (PDB) ({\em https://www.rcsb.org}), since in the OPM database, the mVDAC1 structure has been manipulated exactly like what we need as illustrated in Figure~\ref{boxdomain}, together with  the membrane location numbers $Z1=-12$  \AA \; and $Z2=12$  \AA. This ion channel protein has 4313 atoms and 283 amino acids in one $\alpha$-helix and one $\beta$-barrel (with 19 $\beta$-strands) as depicted in Figure~\ref{structure_3EMN}. The $\beta$-barrel has height 35 \AA  \; and width 40 \AA. The channel pore has the width 27 \AA \; at the entrance and 14 \AA \; at the center.  After downloading a PDB file from the OPM database (with the PDB identification (ID): 3EMN), we converted it to a PQR file on the PDB2PQR web server ({\em http://nbcr-222.ucsd.edu/pdb2pqr\_2.1.1/}) to get the  data missed in the PDB file such as hydrogen atoms,  the atomic charge numbers, and atomic radii. Such a PQR file is a required  input data file for our program package.

\begin{figure}[t]
 \centering
         \begin{minipage}[b]{0.48\textwidth}
                \centering
                \includegraphics[width=\textwidth]{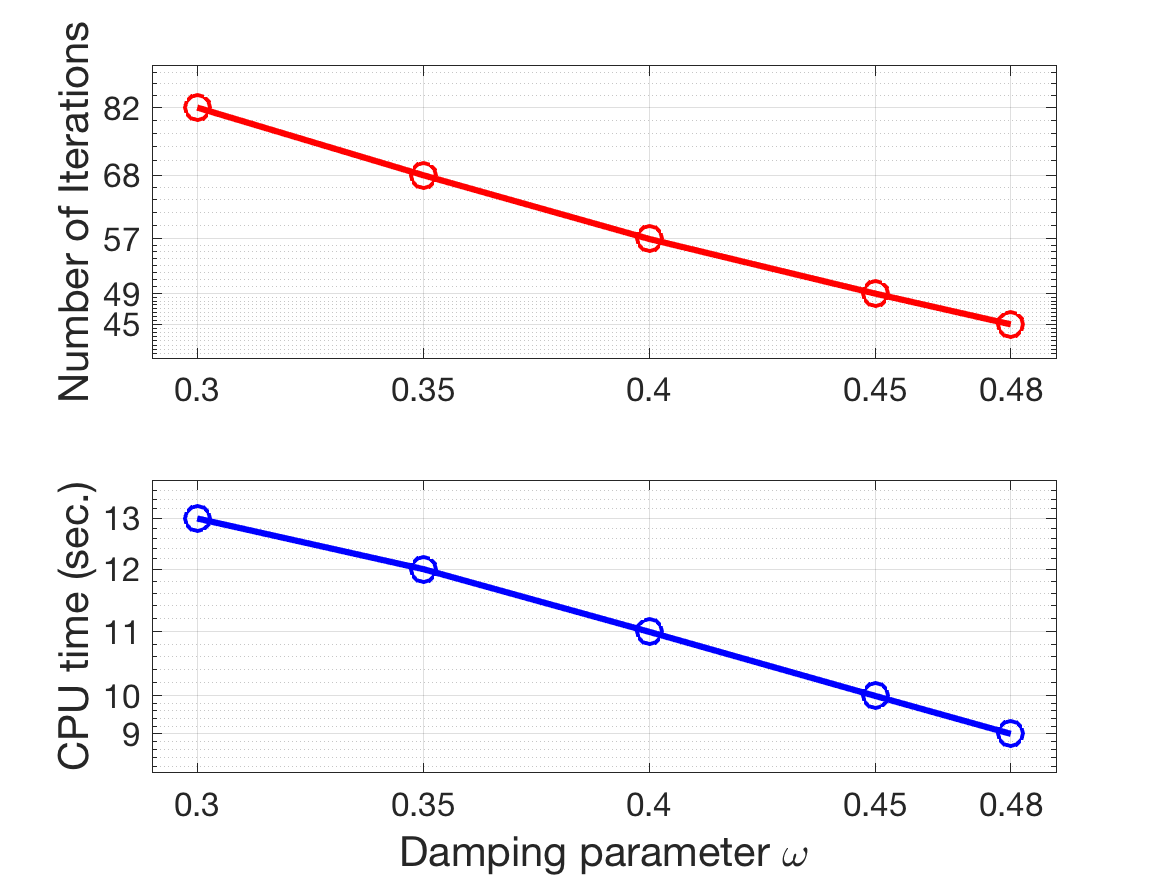}
                \caption{Convergence and performance of the damped two-block iterative method  \eqref{two-block-outloop} for a  finite element system  \eqref{transformed_system} of Model 3 as functions of $\omega$ for a mVDAC1 (PDB ID: 3EMN) in a mixture solution  with  four ionic species Cl$^-$, Na$^+$,  K$^+$, and  NO$_3^-$. }
                \label{converence3emn4ions}
        \end{minipage}
        \quad
        \begin{minipage}[b]{0.48\textwidth}
                \centering
                \includegraphics[width=\textwidth]{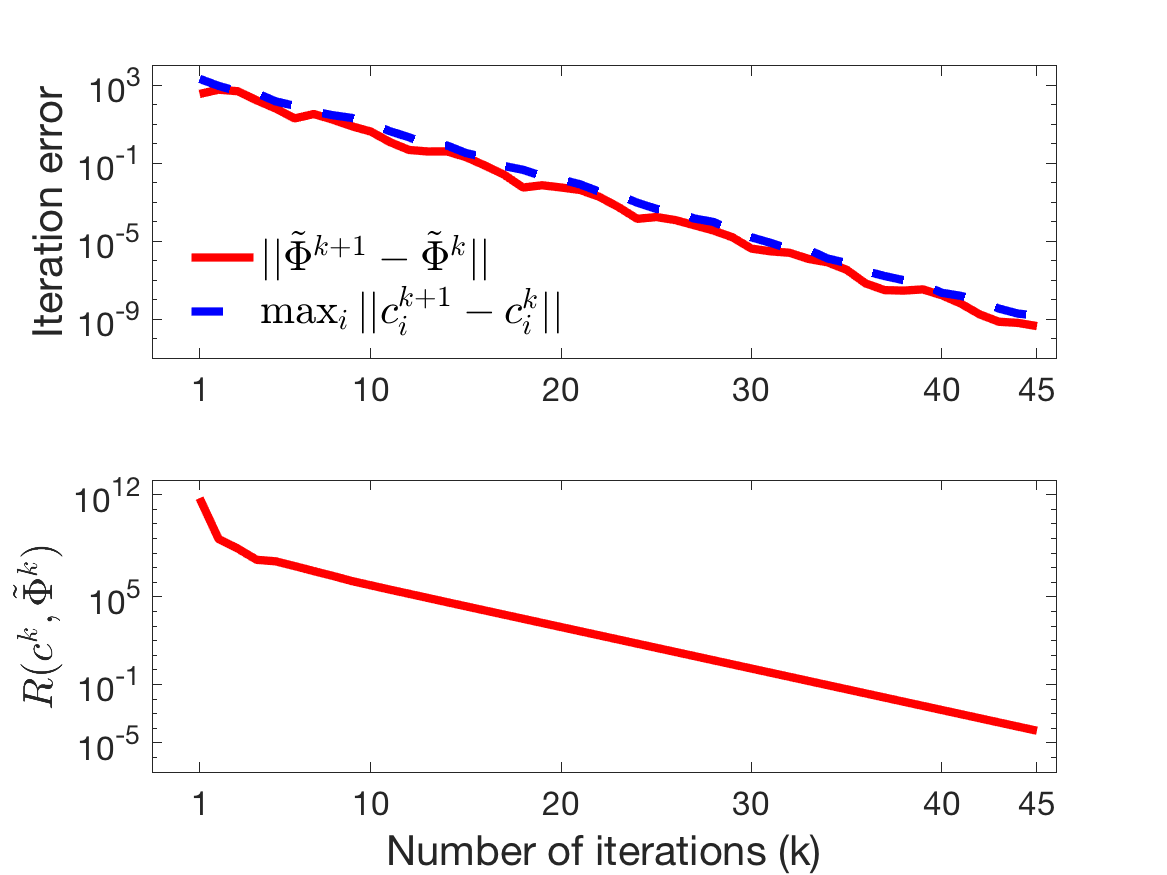}
                \caption{Iteration errors and residual errors  of the damped two-block iterative method  using $\omega=0.48$ for solving a nonlinear  finite element system  \eqref{transformed_system} of Model 3 as functions of the number $k$ of iterations. Here the residual error $R(c,\Phit)$ is defined in \eqref{constrianResidual}. }
                \label{converence3emnNewton}
        \end{minipage}  
\end{figure} 

We constructed a box domain $\Omega$ using $L_{x1} = -36,  L_{x2} = 43,  L_{y1} = -36, L_{y2} = 33,   L_{z1} = -34,$ and $ L_{z2} = 30$ and generated an interface fitted irregular tetrahedral mesh $\Omega_h$ of $\Omega$ (with 49798  mesh points) by an  ion channel finite element mesh program package \cite{ChaoZhenThesis,LiuTT2015}. We then extracted a solvent region mesh $D_{s,h}$ from the box domain mesh $\Omega_h$ (with 29366  mesh points). From Figure~\ref {mesh_3EMN} it can be seen that these two meshes  $\Omega_h$ and $D_{s,h}$ are very irregular due to the complex interfaces $\Gamma_p$ and $\Gamma_{pm}$ or a complex molecular surface of mVDAC1.

\begin{figure}[t]
        \centering
         \begin{subfigure}[b]{0.3\textwidth}
                \centering
                \includegraphics[width=\textwidth]{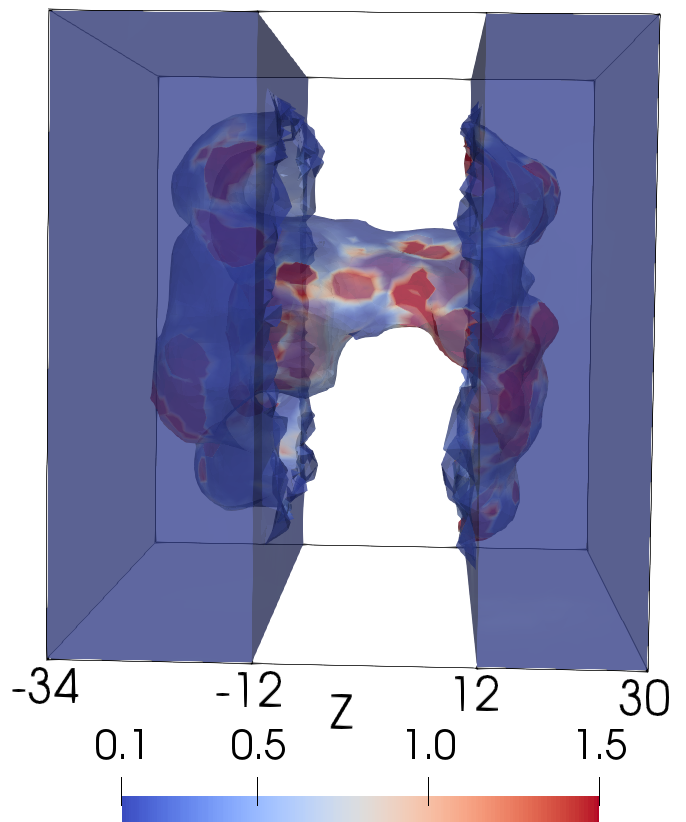}
                \caption{Cl$^-$ concentration}
        \end{subfigure}  
        \qquad \qquad
         \begin{subfigure}[b]{0.3\textwidth}
                \centering
                \includegraphics[width=\textwidth]{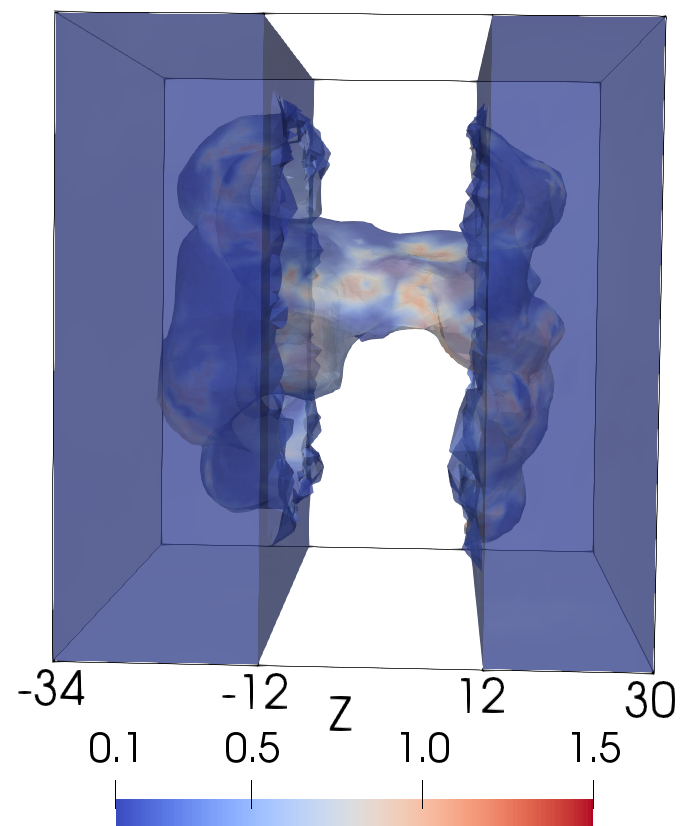}
                \caption{NO$_3^-$ concentration}
        \end{subfigure}
               \begin{subfigure}[b]{0.3\textwidth}
                \centering
                \includegraphics[width=\textwidth]{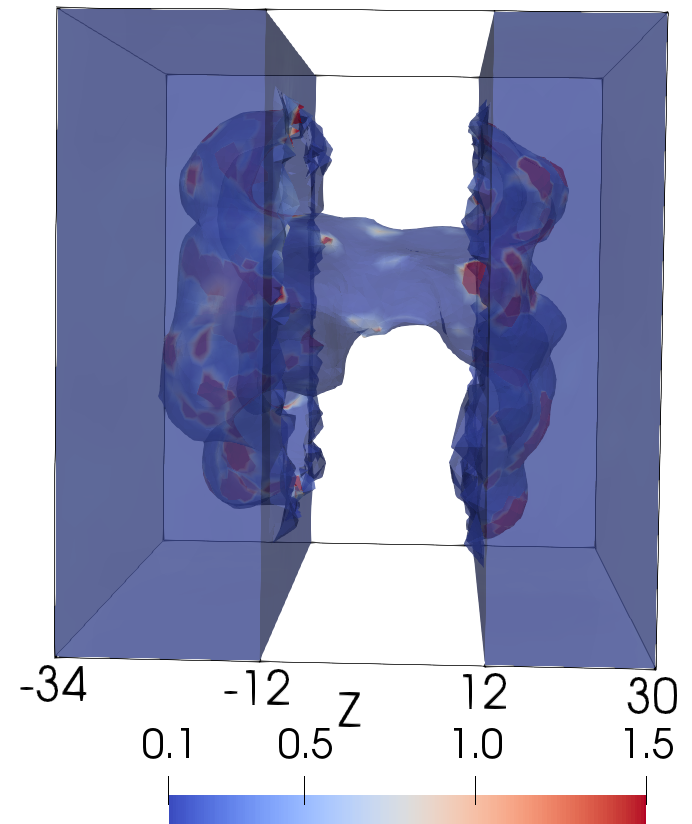}
                \caption{K$^+$ concentration}
        \end{subfigure}  
        \qquad \qquad
         \begin{subfigure}[b]{0.3\textwidth}
                \centering
                \includegraphics[width=\textwidth]{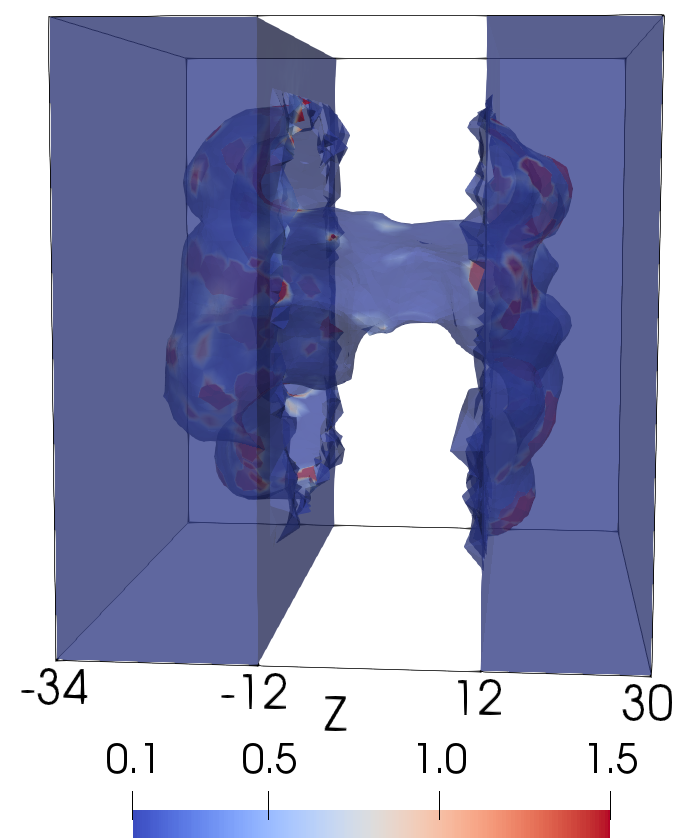}
                \caption{Na$^+$ concentration}
        \end{subfigure}
\caption{The color mappings of the four ionic concentrations generated by the nuSMPBIC finite element program package for mVDAC1 (PDB ID: 3EMN) in a mixture solution of  0.1 molar  KNO$_3$ and 0.1 molar NaCl on one surface view of the solvent region $D_s$ in the $Z$-axis direction.}        
\label{NO3NaClK_test4}          
\end{figure}

\begin{figure}[t]
        \centering
                  \begin{subfigure}[b]{0.45\textwidth}
                \centering
                \includegraphics[width=\textwidth]{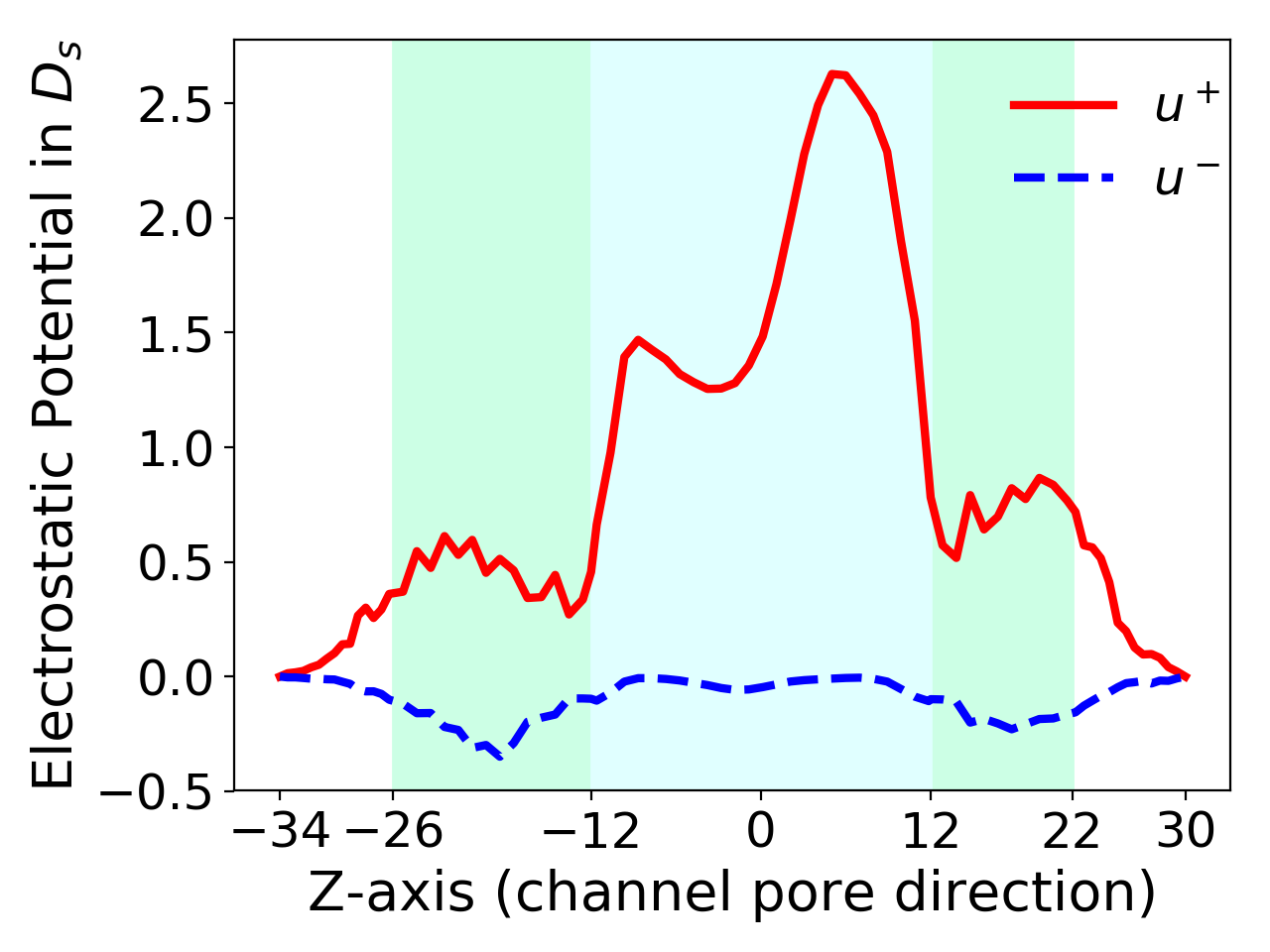}
        \end{subfigure} 
        \qquad
         \begin{subfigure}[b]{0.45\textwidth}
                \centering
                \includegraphics[width=\textwidth]{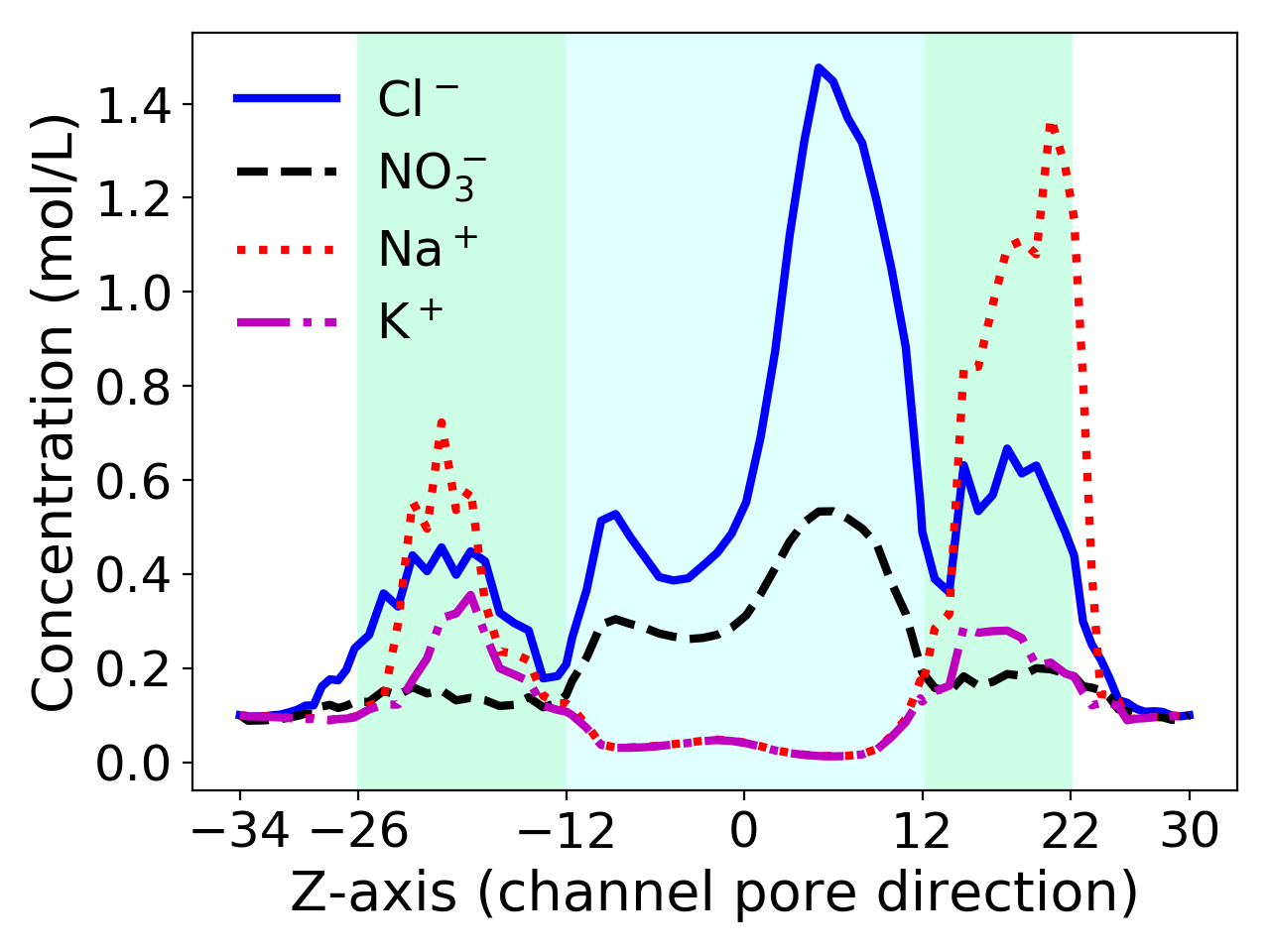}
        \end{subfigure}
\caption{The electrostatic potentials $u^+$ and $u^-$ and ionic concentrations $c_i$ generated by our nuSMPBIC finite element program package for mVDAC1 (PDB ID: 3EMN) in a mixture solution of  0.1 molar  KNO$_3$ and 0.1 molar NaCl  with  four ionic species Cl$^-$, Na$^+$,  K$^+$, and  NO$_3^-$. Here $u^+$ and $u^-$ denote the positive and negative parts of electrostatic potential $u$;  the central part of the channel pore that links to the membrane (i.e., $-12\leq z\leq 12$) is highlighted in light-cyan; and the other parts of the channel pore are highlighted in green.}        
\label{KClNO3Nacase}          
\end{figure}

\begin{figure}
        \centering
         \begin{subfigure}[b]{0.45\textwidth}
                \centering
                \includegraphics[width=\textwidth]{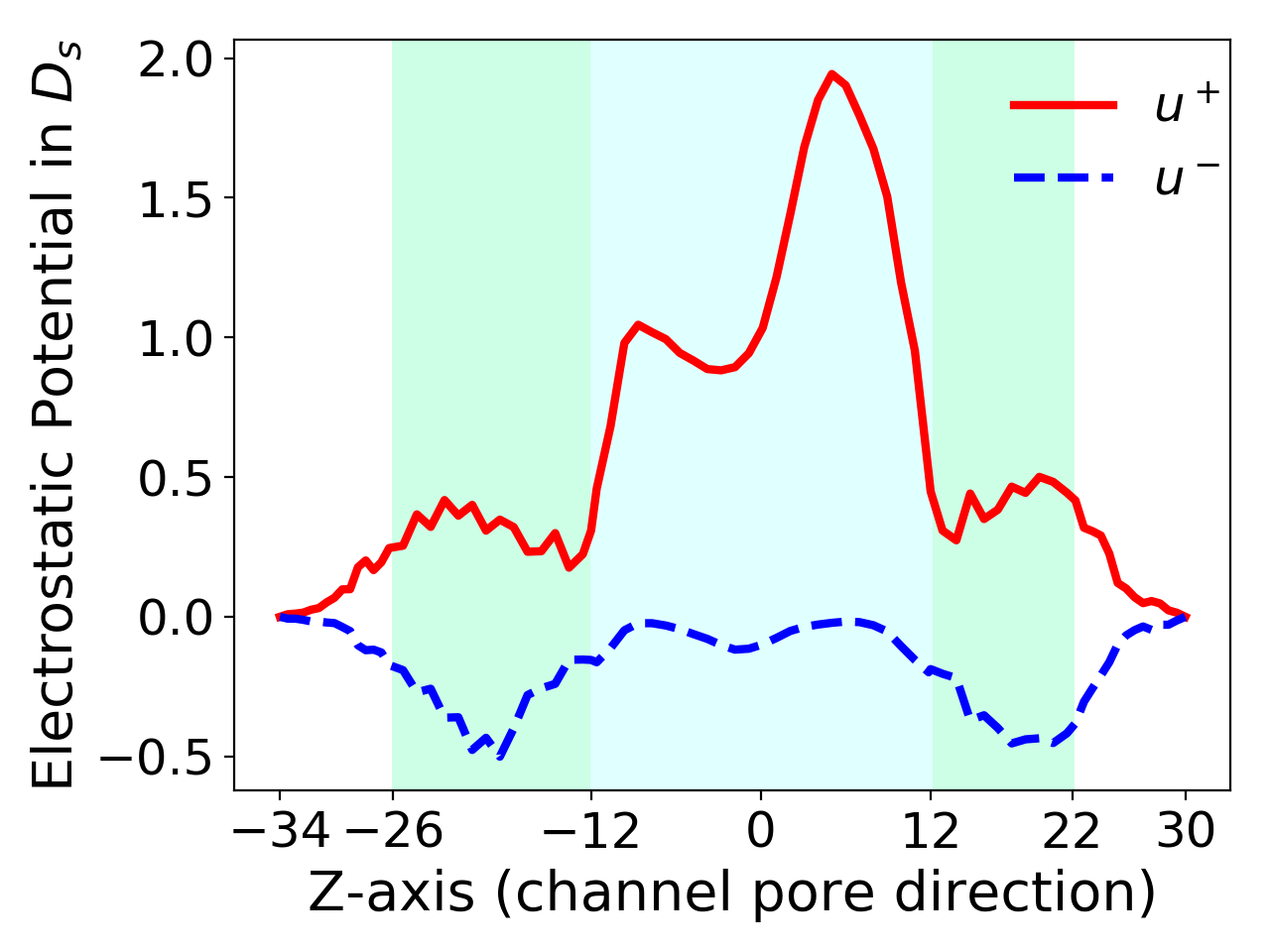}
        \end{subfigure}  
        \qquad
               \begin{subfigure}[b]{0.45\textwidth}
                \centering
                \includegraphics[width=\textwidth]{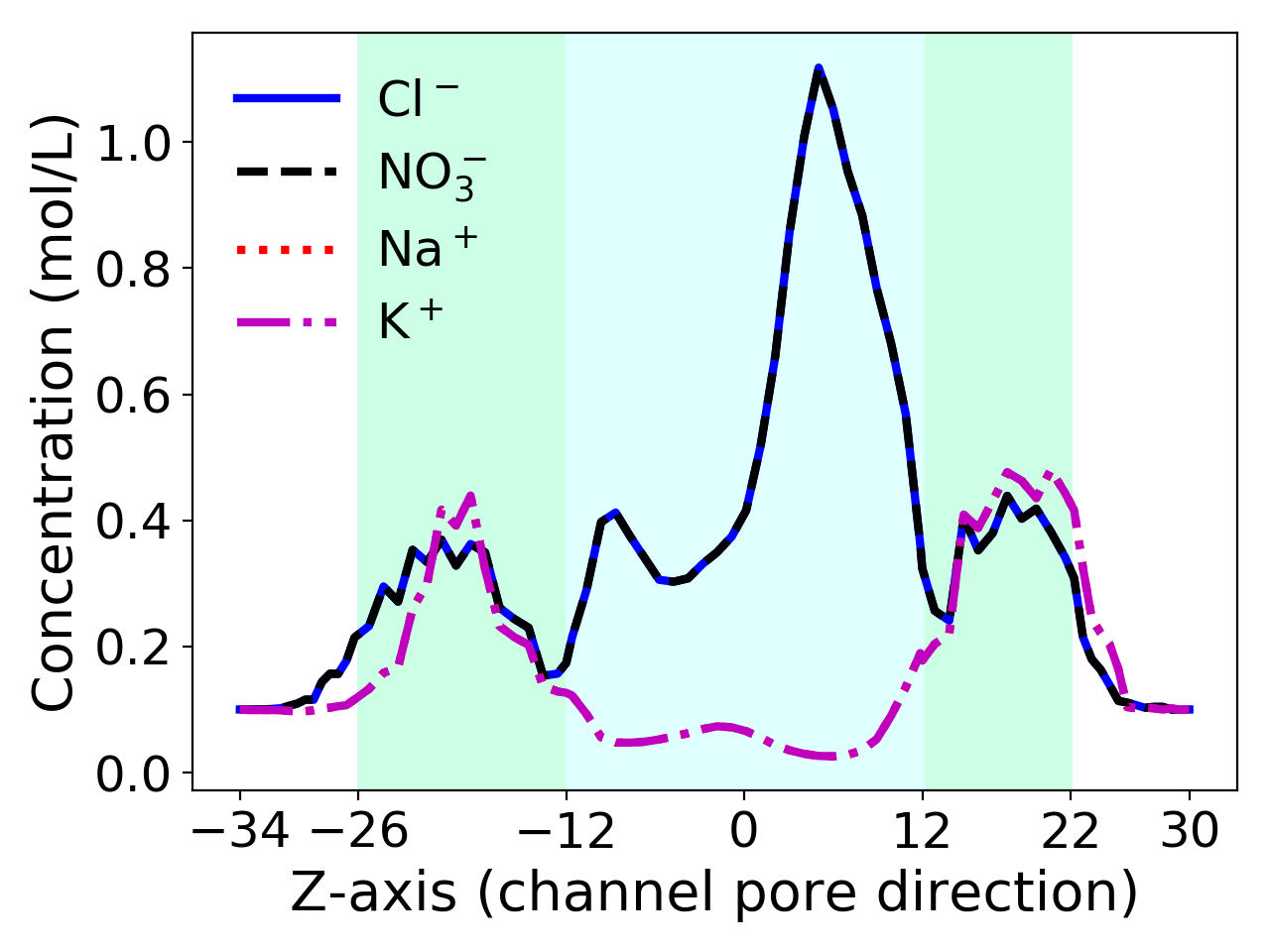}
        \end{subfigure}  
        \caption{The electrostatic potentials $u^+$ and $u^-$ and ionic concentrations $c_i$ generated by the nuSMPBIC finite element program package for mVDAC1 (PDB ID: 3EMN) in a mixture solution of  0.1 molar  KNO$_3$ and 0.1 molar NaCl in a uniform ion size case: all the ions are set to have  the same average volume size $\hat{v}$. That is,  $v_i=\hat{v} = 28.8393$ for $i=1,2,3,4$ for the four ionic species Cl$^-$, NO$_3^-$, K$^+$, and Na$^+$.}        
\label{uSMPBicTests}          
\end{figure}

\begin{figure}[t]
       \begin{subfigure}[b]{0.45\textwidth}
                \centering
                \includegraphics[width=\textwidth]{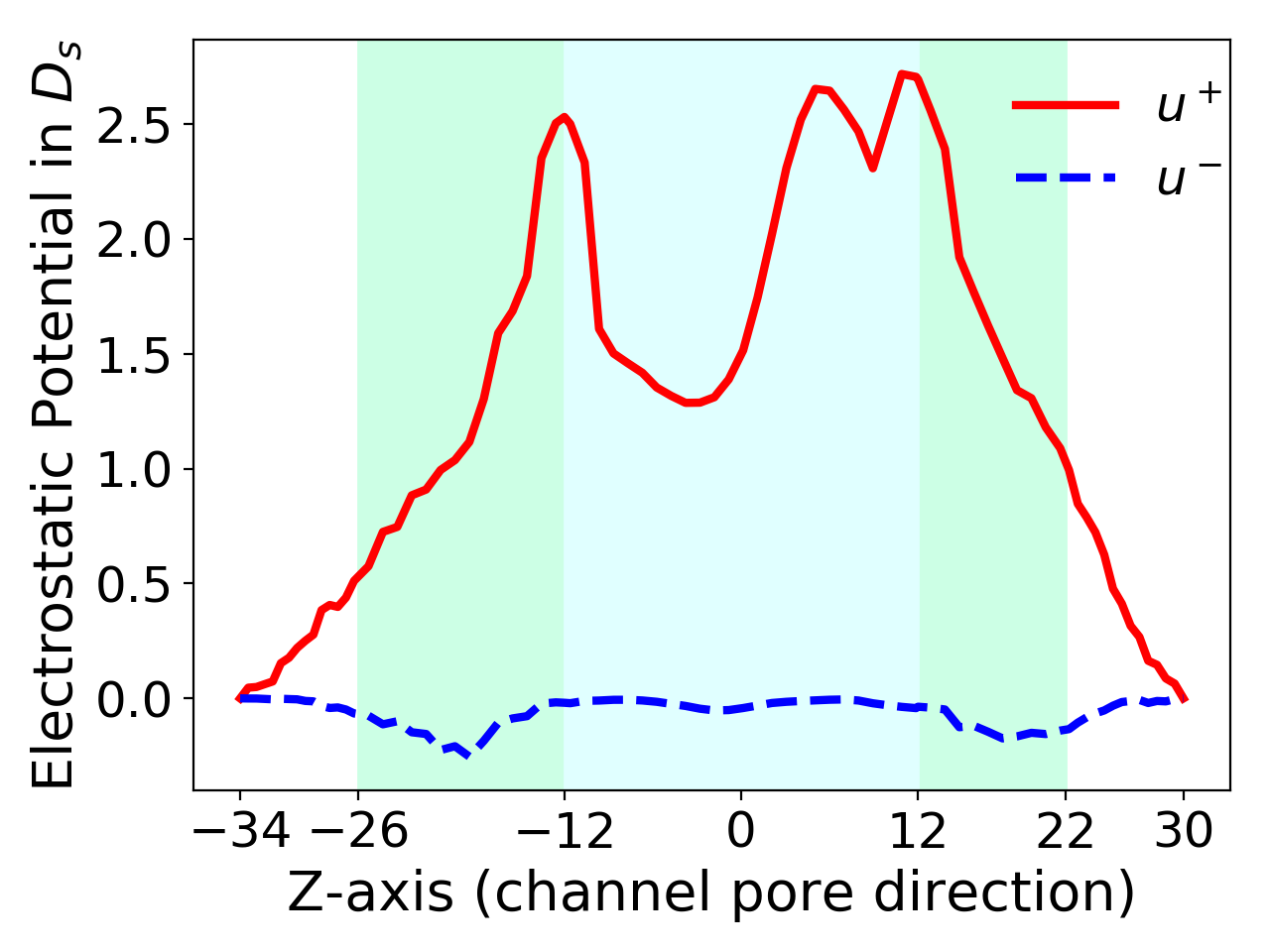}
        \end{subfigure}  
        \qquad
         \begin{subfigure}[b]{0.45\textwidth}
                \centering
                \includegraphics[width=\textwidth]{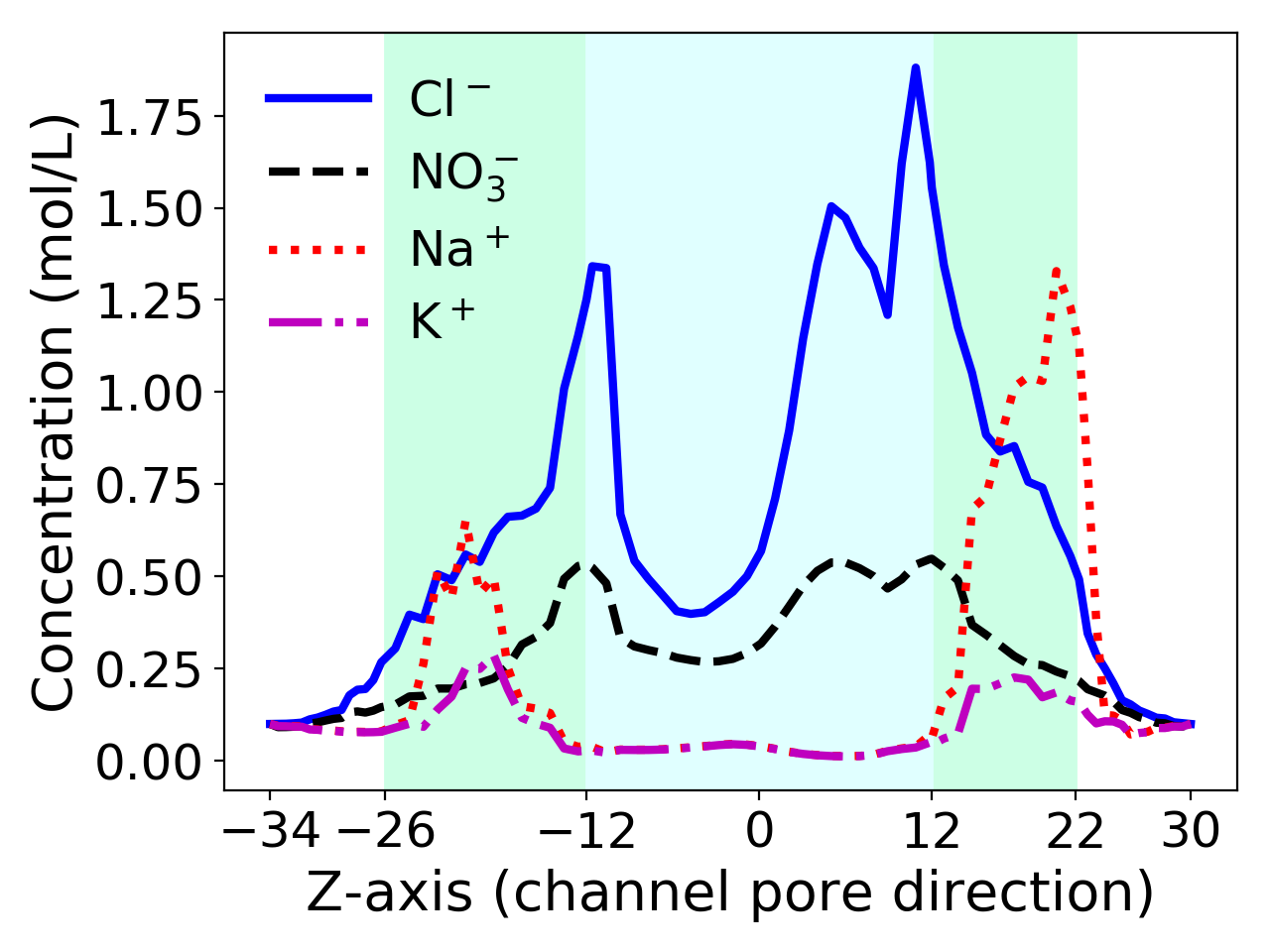}
        \end{subfigure}
\caption{The electrostatic potentials $u^+$ and $u^-$ and ionic concentrations $c_i$ generated by the nuSMPBIC finite element program package with the membrane surface charge density $\sigma=10$ $\mu$C/cm$^2$ for mVDAC1 (PDB ID: 3EMN) in a mixture solution of  0.1 molar  KNO$_3$ and 0.1 molar NaCl.}        
\label{SurfaceChargeCase}          
\end{figure}

All the numerical tests were done on an iMac computer with one 4.2 GHz Intel core i7 processor and 64 GB memory. For simplicity,  we fixed the parameters $\ep = 2$, $\emm=2$,  $\es=80$,  $u_b=0$,  and $u_t=0$ in these tests.
The numerical test results are reported in Table~\ref{table-performance} and Figures~\ref{converence3emn4ions} to \ref{SurfaceChargeCase}.

Figure~\ref{converence3emn4ions} reports the convergence and performance of our damped two-block iterative method \eqref{two-block-outloop} in terms of the number of iterations and computer CPU time. From the figure it can be seen that  the number of iterations and CPU time  were reduced monotonically for  $\omega \leq 0.48$. The damped two-block iterative method was found to be divergent for $\omega > 0.48$. At $\omega=0.48$, it took about 12.63 seconds only to find  one electrostatic potential and four ionic concentrations as a finite element solution of Model 3. These tests  demonstrate the fast convergence of our damped two-block iterative method and the high performance our nuSMPBIC finite element package. 

Figure~\ref{converence3emnNewton} displays a convergence process of the damped two-block iterative method using $\omega=0.48$ in terms two iteration errors $ \|   \Phit^{k+1}  -  \Phit^{k}  \|_{\Omega}$ and   $\max_{1\leq i\leq n} \|   {c}^{k+1}_i -  {c}^{k}_i  \|_{D_s}$ and one residual error $R(c^k,\Phit^k)$ defined in \eqref{constrianResidual}. From the figure it can be seen that this iterative method quickly reduced the two iteration errors from about $10^{3}$ to $10^{-9}$ in 45 iterations and the residual error  from $10^{12}$ to $10^{-5}$, which indicates that the  finite element solution of Model 3 has well reflected the ion size effects since it  satisfied the size constraint conditions \eqref{sizeModifiedBoltzmann} in high accuracy. 

\begin{table}[h]
\caption{A comparison of the performance of our nuSMPBIC finite element solver using the GMRES-ILU iterative method with that using the Gaussian elimination direct method in terms of computer CPU times in seconds. Here $\omega=0.48$.}
\centering 
\begin{tabular}{|c|c|c|c|c|c|}
 \hline
 Linear	  & Calculate  &  Solve Model 2  & Solve \eqref{Phit-iterate0} & Solve Model 3   & Total       \\ 
 solver         &  $G$ \& $\nabla G$   &  for  $\Psi$  & for $\Phit^0$   &  for $(c,\Phit)$ & CPU time  \\    \cline{1-6}
 GMRES-ILU  & 1.13  &  0.91  & 0.60 &   9.37  &  12.01  \\ \hline
Direct  & 1.19 & 2.74 & 2.55  & 12.63 & 19.11   \\ 
 \hline
\end{tabular}
\label{table-performance}
\end{table}

Table~\ref{table-performance} lists the distribution of computer CPU time in the major parts of the nuSMPBIC finite element solver, along with the case using the Gaussian elimination direct method to solve each related linear finite element equation. It shows that the nuSMPBIC finite element solver using the GMRES-ILU iterative method took much less CPU time than that using the Gaussian elimination direct method,  demonstrating the efficiency of the GMRES-ILU method in our nuSMPBIC finite element package. 

Figure~\ref{NO3NaClK_test4} displays the concentrations of ions Cl$^-$, NO$_3^-$, K$^+$, and Na$^+$ on a surface of the solvent region $D_s$ in color mapping. Here all the concentration values more than 1.5 mol/L and less than 0.1 mole/L have been colored in red and blue, respectively.  From Plots (a, b) it can be seen that there are much more anions Cl$^-$ than NO$_3^-$ within the channel pore due to that the ion size of Cl$^-$ is much smaller that that of  NO$_3^-$ (i.e., 24.84 vs. 77.07 in \AA$^3$) even though both Cl$^-$ and NO$_3^-$ have the same charge number $-1$ and bulk concentration $0.1$ mol/L. In the cation case, Plots (c, d) show that there are few cations K$^+$ and Na$^+$ within the channel pore but it is difficult for us to tell which one more or less since each of them only displays one view of  a three-dimensional concentration.

Using a curve visualization scheme reported in \cite{Xie4PNPicNeumann2020}, we calculated the average values $c_i^j$ and $u^{\pm}_j$ of   concentration function $c_i$ and  positive/negative part $u^{\pm}$ of electrostatic potential function $u$  by the formulas
\begin{equation}
\label{values4ci}
   c_i^j = \frac{1}{\|B_{j,h}\|} \int_{B_{j,h}} c_i(\rr) d\rr, \quad u^{\pm}_j = \frac{1}{\|B_{j,h}\|} \int_{B_{j,h}} u^{\pm}(\rr) d\rr,  \quad j=1, 2, \ldots, m, 
\end{equation}
where $m=80$,  $B_{j,h}$ denotes the $j$th block mesh of a solvent region mesh $D_{s,h}$, $u^{\pm} = \frac{1}{2}[u(\rr) \pm | u(\rr) | ]$, and $\|B_{j,h}\|$ denotes the volume of block $B_{j,h}$.
In particular, $B_{j,h}$ was extracted from $D_{s,h}$ by
\[ B_{j,h} = \left( [L_{x1}, L_{x2}] \times [L_{y1}, L_{y2}] \times [z^j -\bar{h} /2,z^j + \bar{h}/2 ] \right)  \cap D_{s,h}, \]
where $z^j$ denotes the $j$th  partition number of the interval $[L_{z1}+ \bar{h}/2, L_{z2} - \bar{h}/2]$ with $\bar{h}=5$ in the $z$-axis direction. These partition numbers include the membrane location numbers $Z_1$ and $Z_2$.  
In addition, we set $z^0= L_{z1}$, $z^{m+1} = L_{z2}$, $ c_i^0 = c_i^b$, $ c_i^{m+1} = c_i^b$, $u^{\pm}_0=u_b$, and $u^{\pm}_{m+1}=u_t$. Using these points $(z^j, c_i^j)$  and  $(z^j, u^{\pm}_j)$ for $j=0, 1,2,\ldots, m, m+1$, we  plotted 2D  curves and displayed them in Figure~\ref{KClNO3Nacase}. Here the central part of the channel pore that links to membrane (i.e., $-12\leq z\leq 12$)  has been highlighted in light-cyan while the other parts of the channel pore in green to let us more clearly view the distribution profiles of $c_i$ and $u^{\pm}$ within the central part of the channel pore.  

From Figure~\ref{KClNO3Nacase}  it can be seen that the positive potential $u^+$ is much stronger than the negative potential $u^-$, causing  anions Cl$^-$ and NO$_3^-$ to  expel most cations K$^+$ and Na$^+$ away from the central portion of a channel pore between $Z1<z<Z2$  highlighted in light-cyan.  We also can see that the maximum concentration value of Cl$^-$ is almost triple that of NO$_3^-$ (1.48 vs. 0.53) simply because the ion size of Cl$^-$  is much smaller than that of NO$_3^-$ (24.84 vs. 77.07). For cations, we find that  the maximum concentration value of Na$^+$ is almost four times that of K$^+$ (1.37 vs. 0.36) in the range $12<z<22$ highlighted in green since  the ion size of K$^+$ is about triple the ion size of Na$^+$ (9.85 vs. 3.59). These test results demonstrate that our nuSMPBIC model can well retain the anion selectivity property of mVDAC1 and that ionic sizes have significant impacts on electrostatics and ionic concentrations. They also indicate that the anion selectivity happens mostly within the central part of the ion channel pore.

The new damped two-block iterative method works for a uniform ion size case too. To confirm it, we did numerical tests using $v_i=\hat{v}$ with $\hat{v} = 28.8393$ and $\omega= 0.35$. The test results were reported in  Figure~\ref{uSMPBicTests}. In this uniform ion size test case, the damped two-block iterative method took 46 iterations and 11 seconds in computer CPU time to satisfy the iteration convergence rule \eqref{Ite-stop}. As shown in Section 6, the nuSMPBIC model can be reduced from a nonlinear system to the nonlinear finite element boundary value  problem \eqref{uSMPhit-def}, which defines the SMPBIC model and has been solved by an efficient modified Newton iterative method in \cite{SMPBEic2019}.
As a comparison, we repeated the test using the SMPBIC finite element package reported in \cite{SMPBEic2019}. In this test, the nonlinear  finite element equation  \eqref{uSMPhit-def} was solved to satisfy the iteration convergence rule \eqref{Ite-stop} in 8 modified Newton iterations but took about 16 seconds in the total CPU time. This comparison test indicates that  the damped two-block iterative method can also be efficient even compared to  the modified Newton iterative method for solving the nonlinear   finite element equation \eqref{uSMPhit-def}. 

However, because all the ions were set to have the same size, the  2D curves of concentrations of the two anion species  Cl$^-$ and NO$_3^-$  and  the  two cation species K$^+$ and Na$^+$ are overlapped each other, respectively,  in Figure~\ref{uSMPBicTests}, implying that the concentrations of  Cl$^-$ and NO$_3^-$ (or  K$^+$ and Na$^+$) are identical each other. This matches exactly what we could expect in physics when all the ions have the same size,  the same charge number ($-1$ for the anions Cl$^-$ and NO$_3^-$ and $+1$ for the cations K$^+$ and Na$^+$),  and the same bulk concentration 0.1 mol/L. Even so, these test results still well reflect  the anion selectivity property of mVDAC1. Hence, the SMPBIC model can remain a valuable model in  ion channel simulation and study due to its simplicity.

Finally, we did tests using the membrane surface charge density  $\sigma=10$ $\mu$C/cm$^2$ to check the affection of membrane charges on a solution of the nuSMPBIC model. The derived potential and concentration functions were reported in Figure~\ref{SurfaceChargeCase}. From a comparison with  those reported in Figure~\ref{KClNO3Nacase} (i.e., the case without considering any membrane charge), we can see that the positive membrane surface charge significantly increased the values of positive electrostatic potential function $u^+$ and anionic concentrations for Cl$^-$ and NO$_3^-$ but had small affection to the negative electrostatic potential function $u^-$ and cationic concentrations for Na$^+$ and K$^+$. This test indicates the importance of considering membrane charges in the calculation of  electrostatic potential and ionic concentration functions.

\section{Conclusions}
We have presented an efficient nuSMPBIC finite element iterative method for solving a nonuniform size modified Poisson-Boltzmann ion channel (nuSMPBIC) model using Neumann-Dirichlet mixed boundary conditions and a membrane surface charge density, along with its finite element program package that works for an ion channel protein with a three-dimensional crystallographic structure, a mixture solution of multiple ionic species, a nonuniform ion size case, and a uniform ion size case. In particular, we divide the nuSMPBIC model into three-submodels, called Models 1, 2, and 3, to overcome the difficulty of solution singularities induced from singular Direct-delta distributions and to sharply reduce the complexity of solving the  nuSMPBIC model. We then have developed an efficient damped two-block iterative method for solving a linear finite element approximation to Model 3 --- a system mixing nonlinear algebraic equations with a finite element variational problem, including a novel modified Newton iterative scheme for solving each related nonlinear algebraic system. Numerical test results on a voltage-dependent anion-channel (VDAC) in a mixture of four ionic species have demonstrated the fast convergence of our damped two-block iterative method, the high performance of our finite element package, and the importance of considering nonuniform ion sizes in the calculation of electrostatic potentials and ionic concentrations. The nuSMPBIC model has also well validated by the anion-selectivity property of VDAC.

Since the focus of this paper is on the presentation of the nuSMPBIC model and its finite element solver, we only reported numerical test results for one ion channel protein in this paper.  To further confirm the effectiveness and performance of our nuSMPBIC finite element solver, we plan to make  numerical experiments on more ion channel proteins and do comparison studies with other ion channel models in our future work. During these new studies, we will further improve the efficiency of the nuSMPBIC finite element iterative method and the quality and usage of the nuSMPBIC finite element package. In this way, our nuSMPBIC finite element package will become a valuable simulation tool not only for quantitative assessment of ion size impacts on ion channel electrostatics and ionic concentrations but also for the study of ion channel selectivity properties, ionic distribution patterns across membrane, and membrane charge effects on electrostatics and ionic concentrations. 

 \section*{Acknowledgements}
 This work was partially supported by  the Simons Foundation, USA,  through research award 711776. 


\end{document}